\documentclass[a4paper,11pt]{article}
\usepackage{amsmath,amstext,amsthm,amssymb,amsfonts,epic,latexsym,hyperref
}
\usepackage[english]{babel}
\usepackage{graphicx}
\usepackage[T1]{fontenc}
\usepackage[latin1]{inputenc}
\usepackage{color}
\usepackage{tikz} 
\usepackage{epsfig}

\newtheorem{mytheo}{Theorem}
\newtheorem{myprop}{Proposition}
\newtheorem{mycorol}{Corollary}
\newtheorem{mylem}{Lemma}

\newtheorem{myremark}{Remark}
\newtheorem{myfact}{Fact}

\newcommand{\nn}{\nonumber}
\newcommand{\8}{\infty}

  {\end{em}\bigskip\par}
  {\bigskip\par}

\newcommand{\fois}{\times}

\renewcommand{\epsilon}{\varepsilon}

\newcommand{\cL}{\mathcal{L}}

\newcommand{\bx}{{\bf x}}

\title{Localization in 
log-gamma polymers with boundaries}

\author{Francis Comets \and  Vu-Lan Nguyen}

\begin{document}

\maketitle

{\footnotesize 
\noindent 
 Universit\'e Paris Diderot--Paris 7, 
Math\'ematiques, 
 case 7012, 75205 Paris
Cedex 13, France.
\\
 	Laboratoire de Probabilit\'es et Mod\`eles Al\'eatoires (LPMA), UMR CNRS 7599
	\\
\noindent e-mail:
\texttt{comets@math.univ-paris-diderot.fr, vlnguyen@math.univ-paris-diderot.fr}

}

\begin{abstract}
Consider the directed polymer in one space dimension in  log-gamma environment with boundary conditions, introduced by Sepp\"al\"ainen~\cite{Sepp12}.
In the equilibrium case, we prove that the end point of the polymer converges in law as the length increases, to a density
proportional to the exponent of a zero-mean random walk. This holds without space normalization, and the mass concentrates in a neighborhood of the minimum of this random walk.
We have analogous results out of equilibrium as well as for the middle point of the polymer with both ends fixed. The existence and the identification of the limit
relies on the analysis of a random walk seen from its infimum.
\\[.3cm]\textbf{Keywords:} directed polymer, random medium, log-gamma model, random walk, localization.
\\[.3cm]\textbf{AMS 2000 subject classifications:}
Primary 60K37, 82D60. Secondary 60K35, 82B41
\end{abstract}

\tableofcontents
\section{Directed polymers and localization}

The directed polymer model was introduced in the statistical physics literature by Huse and Henley~\cite{HuHe85} to mimic the phase boundary of Ising model in presence of random impurities, and it is frequently used to study the roughness statistics of random interfaces. Later, it has been mathematically formulated as a random walk in a random potential by  Imbrie and Spencer~\cite{ImSp88}.  In the $(1+1)$-dimensional lattice polymer case, the random potential is defined by a field of random variables $\{\omega(i,j):(i,j) \in \mathbb{Z}^2\}$ and a polymer $\bx=(x_t; t=0,\ldots n)$ is a nearest neighbor up-right path in $\mathbb{Z}^2$ of length $n$. The weight of a path is equal to the exponent of the sum of the potential it has met on its way. There is 
a competition between the entropy of paths and the disorder strength, i.e., the inhomogeneities of the potential. If the potential is constant,  the path behaves diffusively and spreads smoothly over distances of order of its length. On the contrary, if the potential has large fluctuations,
the path is pinned on sites with large potential values, and it localizes on a few corridors with width of order of unity.  
An early example where this behavior was observed is the parabolic Anderson model yielding a rigourous framework to analyse intermittency~\cite{CaMo94}.   
Recently, significant efforts have been focused on  planar polymer models (i.e. $(1+1)$-dimensional) which fall in the KPZ universality class  (named after Kardar, Parisi and Zhang), see Corwin's recent survey~\cite{CorwinKPZ12}.  In the line of specific first passage percolation models and interacting particle systems, a few explicitly solvable models
were discovered, and they allow for detailed descriptions of new scaling limits and statistics. 
We namely mention  Brownian queues~\cite{MoOc07}, log-gamma polymer~\cite{Sepp12}, KPZ equation~\cite{Hairer13, QuastellKPZ10}. However, the theory of universality classes does not explain the localization phenomena. For instance, the wandering exponent $2/3$ from the KPZ class accounts for the typical transverse displacement of order $n^{2/3}$ of the polymer of length $n$, certainly an important information, however different in nature since it addresses the location of the corridor but not its width.

Let us start by defining the model of directed polymers in random environment.
For each endpoint $(m,n)$ of the path, we can define a point-to-point partition function 
\begin{equation}
\nonumber
 Z_{m,n}^{\omega}=  \sum_{\bx} \exp\Big\{\sum_{t=1}^{m+n}\omega(x_t)\Big\},
\end{equation}
where the sum is over up-right paths $\bx$ that start at $(0,0)$ and end at $(m,n)$. 
The model does not have a temperature in the strict sense of statistical mechanics,
however the parameter $\mu$ entering below the log-gamma distribution of $\omega$ plays a similar role by tuning the strength of the disorder.
The point-to-line partition function is given by 
\begin{equation}
\nonumber
 Z_{n}^{\omega}=  \sum_{k=0}^n Z_{k,n-k}^{\omega}.
\end{equation}
The point-to-line polymer measure of  a path of  length $n$ is 
\begin{equation}
\nonumber
 Q_{n}^{\omega}(\bx)=  \frac{1}{Z_{n}^{\omega}} \exp\Big\{\sum_{t=1}^{n}\omega(x_t)\Big\}.
\end{equation}
It is known that the polymer at a vanishing temperature concentrates on its geodesics. 
However little information is known on the random geodesics~\cite {Newman97},
except under assumptions which are often hard to check \cite{DamronHanson14, GeRASe14}.
A less ambitious way to analyze this localization phenomenon is to consider 
the endpoint of the path, and study the largest probability for ending at a specific point,
\begin{equation}
\label{def:In}
I_n=\max_{x\in\mathbb{Z^d}}Q_{n-1}^{\omega}\{x_n=x\},
\end{equation}
which does not require any information on where the endpoint concentrates.
Observe that $I_n$ is small when the measure is spread out, for example if $\omega$ is constant, but $I_n$ should be much larger when $Q_{n}^{\omega}$ concentrates on a small number of paths. 
In large generality it is proved that the polymer is localized and it is expected from the KPZ scaling that most of the endpoint density lies in a relatively small region around a random point at distance $n^{2/3}$ from the mid-point of the transverse diagonal. The size of this region is much smaller than $n^{2/3}$ and is believed that it is order one.  Moreover, Carmona and Hu~\cite{CaHu02} and Comets, Shiga and Yoshida~\cite{CoShYo03} showed that there is a constant $c_0=c_0(\beta)>0$ such that the event
 \begin{equation*}
\limsup_{n\to\infty} I_n \geq c_0
\end{equation*}
has $\mathbb{P}$-probability one. This property is called
endpoint  localization. In fact, the C\'esaro mean of the sequence $I_n$
is a.s. lower bounded by a positive constant. Analyzing terms in semimartingale decompositions, the technique is quite general, but also very circuitous and thus it only provides rough estimates. 
Recently,  Sepp\"al\"ainen has introduced in~\cite{Sepp12} a new solvable polymer model with a particular choice of the law of the potential. 
In this paper, we consider the log-gamma model, 
taking advantage of its solvability to analyze the mechanism of localization and obtain an explicit 
description. 
The model can be defined with boundary conditions (b.c.), i.e., with a different law for vertices inside the quadrant or on the boundary, see 
(\ref{condition:gamma}).   From now, we will consider this model. First, 
it is convenient to introduce
multiplicative weights
\begin{equation}
\nonumber
Y_{i,j}= e^{\omega(i,j)} , \ (i,j) \in \mathbb{Z}_+^2.
\end{equation}
As discovered in the seminal paper \cite{Sepp12},  some boundary conditions make the model stationary as in Burke's theorem~\cite{OcYor01}, and further, they   make it explicitely sovable.
In this setting, the point-to-point partition function for the paths with fixed endpoint is given by
\begin{equation}
\nonumber
Z_{m,n}= \sum_{\bx \in \Pi_{m,n} }\prod_{t=1}^{m+n} Y_{x_t},
\end{equation}
where $\Pi_{m,n}$ denotes the collection of up-right paths $\bx=(x_t)_{0\leq t\leq m+n}$ in the rectangle $\Lambda_{m,n}=\{0,...,m\}\fois\{0,...,n\}$ that go from $(0,0)$ to $(m,n)$. 
We assign distinct weight distributions on the boundaries $(\mathbb{N}\fois \{0\}) \cup(\{0\}\fois\mathbb{N})$ and in the bulk $\mathbb{N}^2$. In order to make it clear, we use the symbols $U$ and $V$ for the weights on the horizontal and vertical boundaries:
\begin{equation}
\nn
U_{i,0}=Y_{i,0} \  \textrm{and} \  V_{0,j}= Y_{0,j} \ \textrm{for} \  i,j \in \mathbb{N}:=\{1,2,\ldots\}.
\end{equation}
{\bf Model b.c.($\theta)$}: Let $\mu >0$ be fixed. For $\theta \in(0,\mu)$, we will denote by b.c.($\theta$) the model with
\begin{equation}
\label{condition:gamma}
\begin{aligned}
&
\{U_{i,0}, V_{0,j},Y_{i,j}: i,j \in \mathbb{N}\} \  \textrm{are independent with distributions} \\
& U_{i,0}^{-1} \sim \textrm{Gamma}(\theta,1), \quad V_{0,j}^{-1} \sim \textrm{Gamma}(\mu\!-\!\theta,1),\quad Y_{i,j}^{-1} \sim \textrm{Gamma}(\mu,1).
\end{aligned} 
\end{equation}
where $\textrm{Gamma}(\theta,r)$ distribution has density $\Gamma(\theta)^{-1}r^{\theta}x^{\theta-1}e^{-rx}$ with $\theta >0, r>0$.

The polymer model with boundary condition posesses a two-dimensional stationarity property. Using this property, 
Sepp\"al\"ainen~\cite{Sepp12} obtains an explicit expression for  the variance of the partition function, he proves that the  fluctuation exponent of free energy is $1/3$
and that the exponent for transverse displacement of the path is $2/3$.
This model has soon attracted a strong interest: large deviations of the partition function~\cite{GeorSepp13}, explicit formula for the Laplace transform of the partition function at finite size~\cite{CoOcSeppZyg14},  GUE Tracy-Widom fluctuations for $Z_n$ at scale $n^{1/3}$~\cite{BoroCorwReme13}, computations of Busemann functions \cite{GeRASeYi14}.
\medskip

In fact, the model of directed polymers  can be defined in arbitrary dimension $1+d$ and with general environment law, see \cite{ImSp88}, and we now briefly 
mention some results 
for comparison.
In contrast with the above results for $d=1$, if the space dimension is large and the potential has small fluctuations -- the so-called weak disorder regime --,
this exponent is 0, and under $Q_n^{\omega}$ the fluctuation of the polymer path is order $\mathcal{O}(n^{1/2})$ with a Brownian 
scaling limit, see 
 \cite{Bolt89,CoYo06,ImSp88}. More precisely, if the space dimension $d\geq 3$ and if the ratio ${\mathbb E} (e^{2 \omega} ) / ({\mathbb E} e^{\omega})^2$ 
is smaller than the inverse of the return probability for the simple random walk,
%
%
%
%
%
the end point, rescaled by $n^{-1/2}$, converges to 
a centered $d$-dimensional Gaussian vector.
Moreover, under the previous assumptions, $I_n \to 0$ a.s.,   at the rate $n^{-d/2}$ according to the local limit theorem of \cite{Sina95,Varg06}. 
\medskip

Let us come back to the case $d=1$ of up-right polymer paths, more precisely, to the  log-gamma model.
We now give a flavour of our results with an explicit limit description of the endpoint distribution  under the quenched measure.
\begin{equation}
\nonumber
Q_{n}^{\omega}\big\{x_n=(k,n-k)\big\}=\frac{Z_{k,n-k}}{Z_n}, \quad k=0,\ldots, n.
\end{equation}
 For each $n$, denote by  
 \begin{equation}
\label{def:lnor}
 l_n={\arg\max}\{  Z_{k,n-k}; 0 \leq  k \leq n\},
\end{equation}
the location maximizing the above probability, and call it
 the "favourite endpoint". 
\begin{mytheo}\label{theo:main}
Consider the model b.c.($\theta$) with $\theta \in (0,\mu)$.
Define the end-point distribution $\tilde \xi^{(n)}$ centered around its mode, by
\begin{equation} \nonumber
\tilde \xi^{(n)}=( \tilde \xi^{(n)}_k; k \in {\mathbb Z}), \quad \textmd{ with }\  \tilde \xi^{(n)}_k
= Q_{n}^{\omega}\big\{x_n=(l_n+k,n-l_n-k)\big\}.
\end{equation}
Thus, $\tilde \xi^{(n)}$ is a random element of the set ${\mathcal M}_1$ of probability measures on ${\mathbb Z}$. Then, as $n \to \8$, we have convergence in law
\begin{equation}
\label{eq:cvloixi}
 \tilde \xi^{(n)} \stackrel{\mathcal L}{\longrightarrow} \xi \qquad \textmd{ in the space } ({\mathcal M}_1, \| \cdot \|_{TV}),
\end{equation}
where $\|\mu-\nu\|_{TV}= \sum_k |\mu (k)-\nu (k)|$ is the total variation distance.
\end{mytheo}
 The definition of $\xi_k$ is given as a functional of a random walk conditioned to stay positive on 
${\mathbb Z}^+$ and  conditioned to stay strictly positive on ${\mathbb Z}^-$. The explicit
expression for $\xi$ is formula (\ref{eq:thelimit}) below. 
The convergence is not strong but only in distribution.
The above result yields a complete description of the localization phenomenon revealed in~\cite{CaHu02, CoShYo03}.
In particular, the mass of the favourite point is converging in the distributional sense.
\begin{mycorol} \label{cor:first}
Consider the model b.c.($\theta$) from (\ref{condition:gamma}).
With $I_n$ from (\ref{def:In}), it holds
\begin{equation}
\nonumber
I_n \stackrel{\mathcal L}{\longrightarrow} \max\left\{ \frac{\xi_k+\xi_{k+1}}{2}; k \in {\mathbb Z}\right\} >0 ,
\end{equation}
as $n \to \8$. By consequence, $\limsup_n I_n >0  \ \ \mathbb{P}$-a.s.
\end{mycorol}
Moreover, we derive that the endpoint density indeed concentrates in a microscopic region, i.e.,  of size ${\mathcal O}(1)$, around the favourite endpoint.
\begin{mycorol}[
Tightness of polymer endpoint]\label{theo:loc}
Consider the model b.c.($\theta$) from (\ref{condition:gamma}) with $\theta \in (0,\mu)$.
Then we have
\begin{equation} 
\label{eq:endpoint }
\lim\limits_{K\to \infty} \limsup_{n\to \infty} Q_n^{\omega}\big[\| x_n-(l_n,n\!-\!l_n)\|\geq K\big] =0 \ \ 
{\rm in \ probability.}
\end{equation}
\end{mycorol}
Our results call for a few comments.
\begin{myremark} \label{rem:2}
(i) Influence of high peaks in the parabolic Anderson model: it is easy to check that the sequence $Z_{m,n}$ is the unique solution of the parabolic Anderson equation
$$
Z_{m,n}=e^{\omega(m,n)} [Z_{m-1,n}+Z_{m,n-1}]
$$ 
with initial  condition $Z_{0,0}=1$ and  boundary conditions $Z_{-1,n}=Z_{m,-1}=0$. Hence, $Z_{m,n}$ can be interpreted as the mean density at time $m+n$ and location $(m,n)$ of a population starting from one individual at the origin, 
subject to the following discrete dynamics: each particle splits at each integer time into a random number (with mean $2 e^{\omega (m,n)}$ at location $(m,n)$) of identical individual moving independently, and jumping instantaneously 
one step upwards or to the right.  
If $e^{-\omega(m,n)} \sim Gamma(\mu), e^{-\omega(m,0)} \sim Gamma(\theta),$ and
$ e^{-\omega(0,n)} \sim Gamma(\mu-\theta)$, our result applies, and shows that the
 population concentrates around the highest peak and spreads at distance $O(1)$. In particular, the second high peak does not contribute significantly to the measure,
 a feature which is believed to hold in small space dimension only.
 In large time,
 the population density converges, without any scaling, to a limit distribution given by $\xi$.

(ii) Corollary \ref{theo:loc} states uniqueness of the favourite endpoint, in the sense that all the mass is concentrated in the neighborhood of the favourite point $l_n$. 
This property is analogous to uniqueness of geodesics in planar oriented last passage percolation. We refer to \cite{DamronHanson14, GeRASe14} for a detailed and recent account
on this and related questions.
\end{myremark}

Besides the point-to-line polymer measure, we also study in this paper the point-to-point  measure,
for which the polymer endpoint is prescribed. Under this measure,
we obtain similar localization results, that we will state in the next section. They deal with the  location in the direction transverse to the overall displacement
of the "point in the middle" of the polymer chain, and with the "middle edge".
They are the first results of this nature. The main reason is that the general approach via martingales in~\cite{CaHu02, CoShYo03} fails to apply if the 
endpoint of the path is fixed. 
We mention that the alternative method, introduced in \cite{Vargas07} to deal with environment without exponential moments, could be applied to point-to-point measures.
A similar comment holds for another approach, based on integration by parts, which 
has been recently introduced in~\cite{CoCr13} to extend localization results to the path itself -- and then reveal the favourite corridors. So far, it is known to apply to Gaussian environment and Poissonian   environment~\cite{CoYo13}, but whether it covers the log-gamma case is still open.


 As we will see in  section \ref{sec:results}, the localization phenomena around the favourite point in the log-gamma model directly relates to the problem of splitting a random walk at its local minima. This coupling is also the main tool to study the recurrent random walk in random environment \cite{GaPeSh10} in one dimension. In the literature, it was proved by Williams \cite{Will74}, Bertoin~\cite{Bert91b,Bert91a,Bert93}, Bertoin and Doney \cite{BeDo94}, Kersting and Memi\c{s}o\v{g}lu~\cite{KeMe04}  that if the random walk is split at its local minimum, the two new processes will converge in law to certain limits which are related to a process called the random walk conditioned to stay positive/negative. The mechanism is reminiscent of that of the localization in the main valley of the one-dimensional random walk in random environment in the recurrent case,
discovered by Sina\"i \cite{Sinai82} and studied by Golosov~\cite{Golo84}.
\medskip


Our results only hold for boundary conditions ensuring stationarity. A 
 possible way towards the model without boundary conditions
could be via techniques of tropical combinatorics initiated in~\cite{CoOcSeppZyg14}.
\medskip

{\bf Organization of the paper}:
In section \ref{sec:results}, we recall the basic facts on the log-gamma model and state the main localization results both for point-to-line and point-to-point measures. 
 In  section \ref{sec:split}, we introduce the important properties of the random walk conditioned to stay positive that we need to define the limits. In section
 \ref{section: proof} we give
 the proofs of Theorem \ref{theo:main}, Corollaries  \ref{cor:first} and  \ref{theo:loc}.
The last section contains the complete statements for  the point-to-point  measure, together with their proofs.

%

\section{Polymer model with boundary conditions and results} \label{sec:results}

\subsection{Endpoint under the point-to-line measure} \label{subsec:p2l}
Assume the condition (\ref{condition:gamma}). Define for $(m,n)\in\mathbb{Z}_+^2$,
\begin{equation}
\nonumber
U_{m,n}= \frac{Z_{m,n}}{Z_{m-1,n}} \ \textrm{and} \ V_{m,n}= \frac{Z_{m,n}}{Z_{m,n-1}}\; .
\end{equation}

We can associate the $U$'s and $V$'s to  edges of the lattice $\mathbb{Z}_+^2$, so that they represent the 
weight distribution on a horizontal or vertical edge respectively. Let ${\bf e_1}, {\bf e_2}$ denote the unit coordinate vectors in $\mathbb Z^2$. 
For an horizontal edge  $f=\{y-{\bf e_1},y\}$ we set $T_f=U_y$, and  $T_f= V_y$ if $f=\{y-{\bf e_2},y\}$. 
 Let ${\mathbf z}= (z_k)_{k \in \mathbb{Z}}$ be a nearest-neighbor down-right path in $\mathbb{Z}_+^2$, that is, $z_k \in \mathbb{Z}_+^2$ and $z_k-z_{k-1}={\bf e_1} \ \textrm{or} -{\bf e_2}$. Denoting the undirected edges of the path by $f_k=\{z_{k-1},z_k\}$, we then have
\begin{equation*}
T_{f_k}= \left\{ \begin{array}{ll} 
 U_{z_k}, & \textrm{if} \  f_k\   \textrm{is a horizontal edge}\\
  V_{z_{k-1}} & \textrm{if}\   f_k\   \textrm{is a vertical edge.}\\
    \end{array} \right.
\end{equation*}

Sepp\"al\"ainen proved \cite{Sepp12} that the choice of log-gamma distribution provides a stationary structure to the model:
\begin{myfact}[Theorem 3.3 in \cite{Sepp12}]

Assume $(\ref{condition:gamma})$ . For any down-right path $(z_k)_{k\in \mathbb{N}}$ in $\mathbb{Z}_+^2$, the variables $\{T_{f_k}: k\in \mathbb{Z}\}$ are mutually independent with marginal distributions
\begin{equation*}
\begin{aligned}
U^{-1}&\sim \textrm{Gamma}(\theta,1), \      V^{-1}\sim \textrm{Gamma}(\mu\!-\!\theta,1). \end{aligned}
\end{equation*}
\end{myfact}
 By considering the down-right path along the vertices $x$ with $x\cdot ({\mathbf{ e_1+e_2}})=n$, we deduce the following fact, which will be a fundamental ingredient in the next two sections.
\begin{myfact}
 \label{lem:stationary}
For each $n$, the variables $(U_{k,n-k},V_{k,n-k})_{0\leq k\leq n}$ are independent, and
\begin{equation}
\label{result:stationary}
U_{k,n-k}^{-1} \sim \textrm{Gamma}(\theta,1)\ \ \ V_{k,n-k}^{-1} \sim \textrm{Gamma}(\mu\!-\!\theta,1).
\end{equation} 
\end{myfact}


Now, define for each $1\leq k \leq n$ the random variable $X_k^n$ 
\begin{equation*}
X_k^n= -\log (\frac{Z_{k,n-k}}{Z_{k-1,n-k+1}})=-\log (\frac{U_{k,n-k}}{V_{k-1,n-k+1}}),
\end{equation*}
and $X_0^n =0$.
By corollary \ref{lem:stationary}, for each $n$, $(X_k^n)_{1\leq k\leq n}$ are i.i.d random variables, and satisfy
\begin{equation}
\label{condition:favorite_ponit}
\frac{Z_{k,n-k}}{Z_{0,n}}= \exp(-\sum_{i=0}^k X_i^n) .
\end{equation}
Defining $S_k^n= \sum_{i=1}^k X_i^n$,  for $0\leq k\leq n$, we will express the mass at point  $(k,n-k)$  as a function of  $S^n$,
\begin{equation*}
Q_{n}^{\omega}\{x_n=(k,n-k)\}
        = \frac {Z_{k,n-k}}{\sum_{i=0}^n Z_{i,n-i}}
        = \frac{1}{\sum_{i=0}^n \exp(-(S_i^n-S_k^n))}
\end{equation*}

From $(\ref{condition:favorite_ponit})$,  the favourite point $l_n$ defined in (\ref{def:lnor}) is also the minimum of the random walk,
\begin{equation}
\label{params:RW}
l_n={\arg\min}\{  S_k^n; 0 \leq k \leq n\}.
\end{equation}
Since we are only interested in the law of $Q_{n}^{\omega}\{x_n=(k,n-k)\}$, in order to simplify the notion, we consider a single set of i.i.d random variables $(X_k)_{k\in \mathbb{Z}_+}$, with the same distribution under $\mathbb{P}$ as 
$\log(U/V)$,
where 
$U$ and $V$ are independent with the same distribution as in (\ref{result:stationary}).
The associated random walk is given by
\begin{equation} \label{eq:RW}
S_n=\sum_{i=1}^n X_i,
\end{equation}
and we define
\begin{equation*}
\xi_k^n= \frac{1}{\sum_{i=0}^n \exp(-(S_i-S_k))}.
\end{equation*}
Then one can check that for every $n$:
\begin{equation*}
(\xi_k^n)_{0\leq k\leq n} \stackrel{\cL}{=}\big(Q_{n}^{\omega}\{x_n=(k,n-k)\}\big)_{0\leq k\leq n},
\end{equation*}
where $ \stackrel{\cL}{=}$ means equality in  law.
Then instead of considering for each $n$ a new set of i.i.d random variables to calculate $\tilde \xi_k^{(n)}$, we just need the $n$ first steps of the random walk $S_n$ to compute the law of $\xi_k^n$. 
 Hence Theorem 1 can be reformulated as follows:
  \begin{equation}
 \{\xi_{\ell_n+k}^n\}_{k\in\mathbb{Z}}\stackrel{\mathcal{L}}{\longrightarrow}\{\xi_k\}_{k\in \mathbb{Z}} \ , \quad \textmd{ in the } \ell_1-{\rm norm},
 \end{equation}
  with 
  \begin{equation}
\label{eq:l_nencore}
 \ell_n=\underset{k\leq n}{\arg\min}\  S_k.
\end{equation}
Since the environment has a continuous distribution, the minimum is a.s. unique.
The complete construction of the limit $\xi_k$ will be given in Section \ref{section: proof} below in two different cases when $\theta=\mu/2$ and $\theta\ne \mu/2$. 
However, for the convenience of the reader, we give an informal definition, starting with the case $ \theta=\mu/2$. 
Let
 $(S^\uparrow_k,k \geq 0),  (S_k^\downarrow, k \geq 0)$ be two independent
processes, with the first one distributed as the random walk $S$ conditioned to be non-negative (forever), and the second one  
distributed as the random walk $S$ conditioned to be positive (for positive $k$). Since we condition by a negligible event, the proper definition requires some care,
it relies on Doob's $h$-transform. Then, 
 \begin{equation} \label{eq:thelimit}
 \xi_k\qquad = \qquad\left\{ \begin{array}{ll}
 \frac{\displaystyle {\exp(-S_k^\uparrow)}}{\displaystyle 1+\sum_{i=1}^\infty \exp(-S_i^\uparrow)+\sum_{i=1}^\infty \exp(-S_i^\downarrow)}, & \textrm{if} \ k\geq 0,\\
   \frac{\displaystyle \exp(-S_k^\downarrow)}{\displaystyle 1+\sum_{i=1}^\infty \exp(-S_i^\uparrow)+\sum_{i=1}^\infty \exp(-S_i^\downarrow)}, & \textrm{if} \ k<0.  \\
    \end{array} \right.
\end{equation}
In the case $\theta<\mu/2$, then $l_n = {\mathcal O}(1)$, but the limit is still given by the formula \eqref{eq:thelimit}, provided that 
$S^\downarrow$ has a lifetime (equal to the time for the walk to reach its absolute minimum), after which it is infinite. $S^\uparrow$ is  as before,
and it is defined in a classical manner.  Thus, the concatenated process is simply equal to $S$ with a space shift by its minimum value, and time shift by the time 
to reach the minimum.   The last case $\theta>\mu/2$ is similar under the change $k \mapsto n-k$.
\medskip

In particular in the equilibrium case $\theta=\mu/2$, $S_k$ is a random walk with expectation zero. By Donsker's invariance principle,
the random walk has a scaling limit,
\begin{equation} 
 \nonumber
\Big (\frac{1}{\sqrt{n}}S_{[nt]}\Big)_t \stackrel{\mathcal{L}}{\longrightarrow}(W_t)_t\ .
\end{equation}
with $W$ a Brownian motion with diffusion coefficient $2 \Psi_1(\mu/2)$ (there, $\Psi_1(\theta)= (\log \Gamma)"(\theta)$ is the trigamma function). By consequence, the scaling limit of the favourite endpoint is easy to compute in the present model with boundary conditions.
\begin{mytheo}\label{theo:arcsin}
Consider the model b.c.($\theta$) from (\ref{condition:gamma}).\\
(i) When $\theta=\mu/2$, we have
\begin{equation*}
\frac{l_n}{n}\stackrel{\mathcal{L}}{\longrightarrow}\underset{t\in [0,1]}{\arg\min}\  W_t \ ,
\end{equation*}
where 
the limit has the arcsine distribution  with density $\big[\pi\sqrt{s(1-s)}\big]^{-1}$ on the interval $[0,1]$.\\
(ii) When $\theta<\mu/2$, $n-l_n$ converges in law, so
$$\frac{l_n}{n}\stackrel{\mathbb{P}}{\longrightarrow} 1,$$
though when $\theta>\mu/2$, $l_n$ converges in law, so $$\frac{l_n}{n}\stackrel{\mathbb{P}}{\longrightarrow} 0.$$
\end{mytheo} 
In words, the favourite location for the polymer endpoint is random at a macroscopic level  
in the equilibrium case, and degenerate otherwise. Further, the (doubly random) polymer endpoint $x_n$ has the same asymptotics under $Q_n^\omega$, since, by (\ref{eq:endpoint }), $x_n/n$ and $l_n/n$ are
asymptotic as $n \to \8$.  
These results  disagree with KPZ theory, where
the endpoint fluctuates at distance $n^{2/3}$ around the diagonal. A word of explanation is necessary.
The difference comes from the boundary conditions. In the equilibrium case $\mu/2=\theta$
the direction of the endpoint has a maximal dispersion, though in non equilibrium ones it
sticks to one of the coordinate axes. In the model without boundary conditions --that we leave untouched in this paper--, we expect an extra entropy term to come
into the play and balance the random walk $S_n$ in the potential, a factor being of magnitude $n$ and quadratic 
around its minimum (which is the diagonal by symmetry), making the localization happen close to the diagonal and with fluctuations
of order $n^{2/3}$.

\medskip

Finally, we derive a large deviation principle for the endpoint distribution:
\begin{mytheo}\label{theo:deviation}
Consider the model b.c.($\theta$) from (\ref{condition:gamma}).\\
(i) Assume $\theta=\mu/2$. In the Skorohod space $D\Big([0,1],\mathbb{R}^+\Big)$  equipped with Skorohod topology,
 \begin{equation}  \label{equa:deviation}
 \left(\frac{-1}{\sqrt{n}}\log Q_n^{\omega}\big\{x_n=([ns], n\!-\![ns])\big\}\right)_{s\in [0,1]}
 \stackrel{\cL}{\longrightarrow}
 \left(W(s)-\min_{[0,1]} W\right)_{s\in[0,1]}\ .
 \end{equation}
 Moreover, for all segment $A\subset \big\{(s,1-s); s \in [0,1]\big\}$ in the first quadrant,
 \begin{equation}\label{equa:interval}
\frac{-1}{\sqrt{n}}\log Q_n^{\omega}(x_n\in nA)
\stackrel{\cL}{\longrightarrow}
\inf_A W- \min_{[0,1]} W\ .
 \end{equation}
 (ii) Assume $\theta > \mu/2$. Then, as $n \to \8$,
 \begin{equation}\label{equa:deviationb}
- \frac{1}{{n}}\log Q_n^{\omega}\big\{x_n=([ns], n\!-\![ns]) \big\}
 \stackrel{\mathbb{P}}{\longrightarrow}
s \big\vert\Psi_0(\theta)-\Psi_0(\mu-\theta)\big\vert\;,
 \end{equation}
 where $\Psi_0(\theta)= (\log \Gamma)'(\theta)$ is the digamma function.
 Similarly, if $\theta < \mu/2$,
 \begin{equation}\nn
- \frac{1}{{n}}\log Q_n^{\omega}\big\{x_n=([ns], n\!-\![ns])\big\}
 \stackrel{\mathbb{P}}{\longrightarrow}
(1-s) \big\vert\Psi_0(\theta)-\Psi_0(\mu-\theta)\big\vert\;.
 \end{equation}

\end{mytheo}
Then, at logarithmic scale, the large deviation probability for the endpoint is  of order $  \sqrt{n}$
in the equilibrium case, whereas it is of  order $n$ otherwise. This is again specific to boundary conditions, since it is shown in \cite{GeorSepp13} for the model without boundaries that the large deviation probabilities have exponential order $n$ with a rate function which vanishes only on the diagonal ($s=1/2$).
\subsection{Middle point under the point-to-point measure} \label{subsec:p2p}

 In this section, we consider the point-to-point measure with 
 boundary conditions. Fix  $\mu>0$, $(p,q)\in (\mathbb{Z}_+^*)^2$ and for each $N\in \mathbb{N}$, let $R_N$ be the rectangle with vertices $(0,0), (0,qN), (pN,0)$ and $(pN,qN)$. With some fixed
 $$\theta_S, \theta_N \in (0,\mu), \qquad {\rm let} \quad \theta_E=\mu-\theta_N, \; \theta_W=\mu-\theta_S, 
 $$
and denote  $\Theta=(\theta_N, \theta_S, \theta_E, \theta_W)$.
 To sites $(i,j)$ strictly inside $R_N$ we assign  inverse  Gamma variables $Y_{i,j}$ with parameter $\mu$, whereas  to sites on the boundary we assign inverse Gamma variables with parameter $\theta_N, \theta_S, \theta_E$ or $\theta_W $ depending if the boundary is north, south, east or west.

 \medskip
 
 {\bf Model P2P-b.c.($\Theta$)}: Assume
 \begin{equation}
\label{eq:P2Pbc}
\begin{aligned}
&
Y_{i,j}: (i,j) \in R_N\setminus \{{\bf 0}, (pN,qN)\} \  \textrm{are independent with} \\
&  Y_{i,0}^{-1}  \sim \textrm{Gamma}(\theta_S,1) \textmd{ for } i \! \in \![1, pN], \quad  Y_{pN,j}^{-1} \sim \textrm{Gamma}(\theta_E,1) \textmd{ for } j \!\in \![1, qN\!-\!1],\\
&  Y_{0,j}^{-1}  \sim \textrm{Gamma}(\theta_W,1) \textmd{ for } j \!\in  \! [1, qN], \quad \! \! Y_{i,qN}^{-1} \sim \textrm{Gamma}(\theta_N,1) \textmd{ for } i \! \in \![1, pN\!-\!1],\\
& Y_{i,j}^{-1} \sim \textrm{Gamma}(\mu,1) \textmd{ for }
1\leq i \leq pN-1 \textmd{ and } 1\leq j \leq qN-1.
\end{aligned} 
\end{equation}

The point-to-point polymer measure is the probability measure on $\Pi_{pN,qN}$ given by
\begin{equation}
\nn
Q_{pN,qN}^{\omega}(\bx)=   \frac{1}{Z_{pN,qN}^{\omega}} \exp\Big\{\sum_{t=1}^{(p+q)N-1}\omega(x_t)\Big\}.
\end{equation}
For a path $\bx \in \Pi_{pN,qN}$ denote by 
$$t^-=\max\{t: x_t \cdot (q{\bf e_1}+p{\bf e_2}) \leq pqN\}$$ the "time it crosses the transverse diagonal". The 
coordinate of the crossing point can be described up to a multiplicative factor by the integer
\begin{equation}
\label{eq:transv}
F(\bx)= (x_{t^-}+x_{t^-+1})\cdot  (q{\bf e_1}-p{\bf e_2}) .
\end{equation}

\begin{mytheo}\label{theo:introP2P}
Consider the model P2P-b.c.($\Theta$). Then, there exist a random integer $m_N$ depending on $\omega$ and a random probability measure $\hat \xi$ on $\mathbb Z$ such that, as $N \to \8$,
\begin{equation}  \nonumber
\Big( Q_{pN,qN}^{\omega}( F(\bx)= m_N + k); k \in {\mathbb Z} \Big)
  \stackrel{\mathcal L}{\longrightarrow} \hat \xi,
\end{equation}
 in the space $({\mathcal M}_1, \| \cdot \|_{TV})$.
\end{mytheo}
We recall that middle-point localization for the point-to-point measure is
 not covered by the usual martingale approach to localization, and this result is the first one of this nature.
 Here also the limit can be described in terms of the minimum of a functional of random walks. The appropriate form of the claim and the limit itself are given
 in Theorem \ref{th:last}, section \ref{sec:last}. We end with a complement.

\begin{mytheo}\label{theo:arcsin2}
Consider the model P2P-b.c.($\Theta$), and recall $m_N$ from Theorem \ref{theo:introP2P}. \\

(i) When $\theta_N=\theta_S$,  as $N \to \8$, 
\begin{equation}  \label{eq:s3}
\frac{m_N}{4Npq}+\frac{1}{2}\stackrel{\mathcal{L}}{\longrightarrow} \underset{t\in [0,1]}{\arg\min}\  W_t \ ,
\end{equation}
where 
the limit has the arcsine distribution. \\

(ii) When $\theta_N < \theta_S$, then $m_N+2pqN$ converges in law, so
$$\frac{m_N}{4Npq}+\frac{1}{2}
\stackrel{\mathbb{P}}{\longrightarrow} 0,$$
but when $\theta_N > \theta_S$, $m_N-2pqN$ converges in law, so
 $$\frac{m_N}{4Npq}+\frac{1}{2}
\stackrel{\mathbb{P}}{\longrightarrow} 1,$$
\end{mytheo} 
We stress that the equilibrium relation (\ref{eq:s3}) holds whatever $p$ and $q$ are, provided that $\theta_N=\theta_S$.
\bigskip

In order to prove all these results, the direct approach is to understand the growth of the random walk  seen from its local minima. In the next section, we will present different results about the decomposition of random walk around its minima.


\section{Splitting at the infimum 
 and random walk conditioned to stay positive } \label{sec:split}
  Through this section, we will only consider the equilibrium case $\theta=\mu/2$, i.e. when the random walk $S=(S_k, k\geq 0)$ in (\ref{eq:RW}) has mean 0. 
The problem of path decomposition for Markov chains at its infimum points is well studied in the literature by  Williams~\cite{Will74}, Bertoin~\cite{Bert91b, Bert91a, Bert93, BeDo94}, Kersting and Memi\c{s}o\v{g}lu~\cite{KeMe04}. We will follow the fine approach of Bertoin~\cite{Bert93}. We mention at this point that the case of a walk drifting to infinity was considered by Doney~\cite{Doney89}. However for our purpose, we do not need such sophisticated results when $\theta$ is different from $\mu/2$.

 First we will introduce the random walk conditioned to stay non negative and explain how it relates to the decomposition of random walk at its minimum. Then we present Tanaka's construction and its consequence on the growth of the walk around the minimum.
 
 \subsection{Random walk conditioned to stay non negative}
Recall that $S_0=0$.  Define the event that the random walk stays non negative $$\Lambda=\{S_k\geq 0\  \textmd{for all}\ k\geq 0\}.$$
  As the random walk does not drift to $+ \infty$ this event has probability  $\mathbb{P}[\Lambda]=0$. In order to give a meaning for the conditioning with respect to $\Lambda$, we can approximate $\Lambda$ with some other event  $\Lambda_n$. The natural choice here is 
  \begin{equation}  \nonumber
  \Lambda_n=\{S_k\geq 0, \forall\ 0\leq k \leq n\}.
  \end{equation} 
 and we would like to study the  asymptotics for large $n$ of the law of $S$ conditioned by $\Lambda_n$.
 
 Let us  introduce some basic notation. For every real number $x$, we denote by $\mathbb{P}_x$ the law of the random walk $S$ started at $x$, and we put $\mathbb{P}=\mathbb{P}_0$. Let $\tau$ be the first entrance time in $(-\infty,0)$:
 $$\tau=\min\{k\geq 1:S_k<0\}.$$
 In particular $\Lambda_n=\{\tau>n\}$. Let $(H,T)=((H_k,T_k),k\geq 0)$ be the strict ascending ladder point process of the reflected random walk $-S$. That is, $T_0=0$ and, for $k=0,1,\ldots$,
 $$H_k=-S_{T_k}, \ \ T_{k+1}=\min\{j>T_k:-S_j>H_k\},$$
 with the convention $H_k\!=\!\infty$ when $T_k\!=\!\infty$. The variable $H_1$ is called the~first strict ascending height of $-S$, they are depicted in  Figure 
 1. The renewal function associated with $H_1$ is
 $$
 V(x)=\sum_{k=0}^{\infty} \mathbb{P}(H_k\leq x).
$$
\begin{center} 
	\begin{figure}\label{Fig:ladder}
		\begin{tikzpicture}
		  		\draw[->,black, ultra thick] (0,0)--(11,0);
		  		\draw[->,black, ultra thick] (0,0)--(0,3);
		  		
		  		\draw[black, ultra thick] (0,0)--(0.5,-0.3)--(1,-1.2)--(1.5,-0.5)--(2,-0.7)--(3,-0.3)--(3.5,0.5)--(4,0.2)--(4.5,-0.9)--(5,-1.2)--(5.5,-1.5)--(6,-1.3)--(6.5,-1.5)--(7,-2)--(7.5,-1.6)--(8,-0.7)--(8.5,-0.5)--(9,-0.3)--(9.5,0.2)--(10,1.3)--(10.5,2.2)--(11,2);
		  		
		  		\draw [dashed] (0, 0.5) -- (3.5,0.5);
				\draw [dashed] (0, 1.3) -- (10,1.3);
		  		\draw [dashed] (0, 2.2) -- (10.5,2.2);
		  		\draw [dashed] (3.5, 0) -- (3.5,0.5);
		  		\draw [dashed] (3.5, 0) -- (3.5,0.5);
		  		\draw [dashed] (10, 0) -- (10,1.3);
		  		\draw [dashed] (10.5, 0) -- (10.5,2.2);
		  		
		  		\node at (3.5,-0.5){\small{{\bf $T_1$}}};
		  		\node at (10,-0.5){\small{{\bf $T_2$}}};
		  		\node at (10.5,-0.5){\small{{\bf $T_3$}}};
		  		\node at (-0.5,0.5){\small{{\bf $-H_1$}}};
				\node at (-0.5,1.3){\small{{\bf $-H_2$}}};
				\node at (-0.5,2.2){\small{{\bf $-H_3$}}};
		\end{tikzpicture}
			\caption{ \small  The strict ascending ladder of the random walk $-S$. Line segments represent jumps. $T_1$ is the first time the walk is positive, 
			$-H_1=-S_{T_1}$ is the value. $T_2$ is the next time the walk takes a larger value, denoted by $-H_2$. Etc.\dots}
	\end{figure}
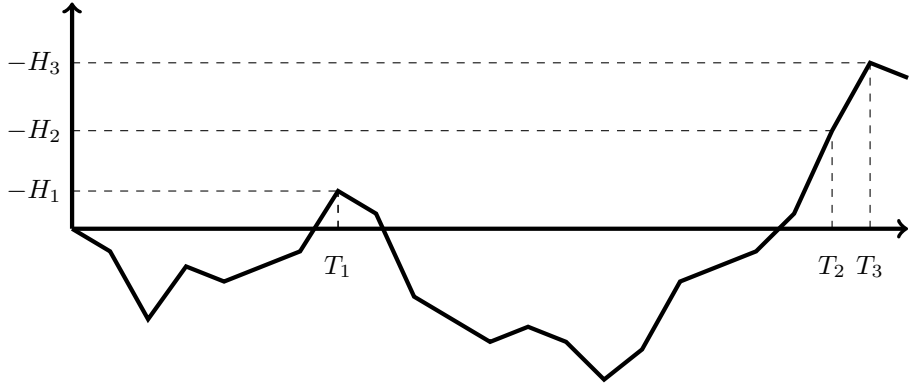
\end{center}
 By the duality lemma (\cite{Feller}, Sect. XII.2), we can rewrite the renewal function for $x \geq 0$ as
 \begin{equation}
 \label{eq:V}
V(x)= 1+ \mathbb{E}(\sum_{i=1}^{\sigma(0)-1}1_{\{-x\leq S_i\}}),
\end{equation}
 where 
 \begin{equation}  \nonumber
 \sigma(0)=\min\{k\geq 1: S_k\geq 0\}.
\end{equation}
 Now we define  Doob's h-transform $P_x^V$ of $P_x$ by the function $V$, i.e.,  the law of the homogeneous Markov chain on the nonnegative real numbers with transition function:
 \begin{equation} \label{eq:Doob}
p^V(x,y)=\frac{V(y)}{V(x)}p(x,y)1_{\{y \geq 0\}}.
\end{equation}
 Here $p, P_x$  denote the transition density and the law of the chain $S$ starting from $x$.
 By definition,  if $f(S)=f(S_0,S_1,...,S_k)$ is a functional depending only on the $k$ first steps of the random walk, then
 $$\mathbb{E}_x^V[f(S)]=\frac{1}{V(x)}\mathbb{E}_x[V(S_k)f(S),k<\tau].$$ 
 (We use the standard notation $E[Z,A]=E[Z {\mathbf 1}_A]$ for an integrable r.v. $Z$ and an event $A$.)
 We denote by $(S_k^{\uparrow})_{k \geq 0}$ the chain starting from 0, 
 \begin{equation}
\label{eq:Sup}
\mathbb{E} (f(S_1^{\uparrow},\ldots,S_k^{\uparrow}))=
 \mathbb{E}_0^V(f(S)).
\end{equation}
 The following result shows that it yields the correct description of the random walk conditioned to stay non negative.
 \begin{myprop}
 \label{theo:convergence2 }
 For a bounded Borel function $f(S)=f(S_1,\ldots,S_k)$,
 \begin{equation}  \nonumber
 \lim_{n\to\infty} \mathbb{E}(f(S)\vert \Lambda_n)=\mathbb{E}(f(S^{\uparrow})).
 \end{equation}
 \end{myprop}
 
  \textit{Proof.} First we will prove the following lemma: 
 \begin{mylem}
 \label{lem:ratio}
 For every $x\geq 0$, we have 
 \begin{equation}  \nonumber
 \liminf_{n\to \infty} \frac{\mathbb{P}_x(\Lambda_n)}{\mathbb{P}(\Lambda_n)}\geq V(x).
 \end{equation}
 \end{mylem}
 \textit{Proof of Lemma \ref{lem:ratio}.} 
 Recall that $(H_k,T_k)$ denotes the $k$th ascending ladder point of $-S$. Let $a_{n,k}$ be the event  $\{H_k\leq x,T_k\leq n, T_{k+1}-T_k>n\}$. On the event $a_{n,k}$, we have  $\max_{k\in [0,n]}\{-S_k\}=H_k\leq x$. It implies that $\min_{k\in [0,n]}S_k+x\geq 0$, and by consequence $\tau>n$ under $\mathbb{P}_x$. Moreover the events $a_{n,k}$ are clearly disjoint, then we have:
 $$\mathbb{P}_x(\tau>n)\geq \sum_{k=0}^{\infty} \mathbb{P}(a_{n,k}).$$
 By the Markov property at $T_k$, we have:
 $$\mathbb{P}_x(\tau>n)\geq \mathbb{P}(\tau>n)\sum_{k=0}^{\infty} \mathbb{P}(H_k\leq x, T_k\leq n).$$
 By monotone convergence,
 \[
 \lim_{n\to \infty} \sum_{k=0}^{\infty} \mathbb{P}(H_k\leq x, T_k\leq n)=\sum_{k=0}^{\infty}\mathbb{P} (H_k\leq x)=V(x).
 \]
 which yields the lemma. \qed
 \medskip

 Now, we can complete the proof of Proposition \ref{theo:convergence2 }.
 Without loss of generality we may assume that $0\leq f \leq 1$. By the Markov property, for $k\leq n$,
 $$\mathbb{E}(f(S),\Lambda_n)=\mathbb{E}\big(f(S)\mathbb{P}_{S_k}(\Lambda_{n-k}),\tau>k\big)\geq \mathbb{E}\big(f(S)\mathbb{P}_{S_k}(\Lambda_{n}),\tau>k\big).$$
 We deduce from Lemma \ref{lem:ratio} and Fatou's lemma that
 $$\liminf_{n\to \infty} \mathbb{E}(f(S)\vert \Lambda_n)\geq \mathbb{E} (f(S)V(S_k),\tau>k)= \mathbb{E}_0^V(f(S))$$
 since $V(0)=1$.
 Replacing $f$ by $1-f$, we get
 \begin{eqnarray*}
 \limsup_{n\to \infty} \mathbb{E}(f(S)\vert \Lambda_n) &=&1-\liminf_{n\to \infty} \mathbb{E}((1-f)(S)\vert \Lambda_n)
 \\
 &\leq & 1-\mathbb{E}_0^V((1-f)(S))= \mathbb{E}_0^V(f(S)),
\end{eqnarray*}
 which completes the proof of Proposition \ref{theo:convergence2 }.
 \qed
 \medskip

 Now we will show that the random walk conditioned stay non negative is the natural limit for the random walk seen from its local minima. Recall $\ell_n$ from (\ref{eq:l_nencore}). The following property is crucial.
 
 \begin{myprop}
 \label{theo:finite}
 For a bounded function $f(x_1,\ldots ,x_k)$ and $\epsilon\in(0,1)$, we have
 \begin{equation} 
  \nonumber
 \lim_{n\to\infty} \mathbb{E}\big(f(S_{\ell_n+1}\!-\!S_{\ell_n},
\ldots,S_{\ell_n+k}\!-\!S_{\ell_n})\vert n-\ell_n>n\epsilon\big)= \mathbb{E}_0^V(f(S)).
 \end{equation}
 \end{myprop}
 \textit{Proof.} 
 We have
 \begin{eqnarray}\nonumber
 \mathbb{E}\big[f(S_{\ell_n+1}\!-\!S_{\ell_n},\ldots,S_{\ell_n+k}\!-\!S_{\ell_n}), {n-\ell_n>n\epsilon}\big]
 \qquad \qquad \qquad \\  \label{eq:rev1}
 \qquad \qquad = \sum_{i=0}^{ \lceil n\!-\!n\epsilon\rceil-1}\mathbb{E}\big[f(S_{i+1}\!-\!S_{i},\ldots,S_{i+k}\!-\!S_{i}), {\ell_n=i}\big].
 \end{eqnarray}
 
 On the other hand we can write the event $\{\ell_n=i\}$ as
 $$\{\ell_n=i\}=\{S_j\geq S_i, \forall j\leq i\}\cap \{S_j\geq S_i, \forall  j\in [ i,n]\}.$$
 Both random variables $f(S_{i+1}-S_{i},\ldots,S_{i+k}-S_{i})$ and $1_{\{S_j\geq S_i, \forall  j\in [ i,n]\}}$ are measurable 
with respect  to $\sigma(X_{i+1},\ldots,X_n)$ and thus are independent of the event $\{S_j\geq S_i, \forall j\leq i\}$ which is in $\sigma(X_1,\ldots,X_i)$. Then we obtain:
 \begin{eqnarray*}
&\mathbb{E}\big[f(S_{i+1}-S_{i},\ldots,S_{i+k}-S_{i}), {\ell_n=i}\big]=  \qquad \qquad \qquad \qquad \qquad \qquad \qquad \qquad \\ & \qquad \qquad 
\mathbb{E}\big(f(S_{i+1}-S_{i},\ldots,S_{i+k}-S_{i})1_{\{S_j\geq S_i, \forall  j \in  [i,n]\}}\big) 
\times \mathbb{P}(S_j\geq S_i, \forall j\leq i).
\end{eqnarray*}
 Applying the Markov property at time $i (0<i<n-n\epsilon)$, we obtain
  \begin{eqnarray}\nonumber
\mathbb{E}\big[ f(S_{i+1}-S_{i},...,S_{i+k}-S_{i}), {\ell_n=i} \big] \qquad \qquad \qquad \qquad \\
\label{eq:rev2}  \qquad \qquad
=\mathbb{E}\big[ f(S_{1},\ldots,S_{k}), {\tau> n-i}\big]\times\mathbb{P}(S_j\geq S_i, \forall j\leq i).
\end{eqnarray}
 From Proposition \ref{theo:convergence2 }, for fixed $\delta>0$, there exists $n(\delta)$ such that for $n\geq n(\delta)$, 
   \begin{eqnarray} 
\big \vert \mathbb{E}\big[f(S_{1},\ldots,S_{k}) \vert {\tau> n\epsilon}\big]- \mathbb{E}_0^V[f(S_1,\ldots,S_k)]\big\vert \leq \delta. \label{eq:rev3}
\end{eqnarray}
Combining \eqref{eq:rev1}, \eqref{eq:rev2} and \eqref{eq:rev3}, we obtain
 \begin{eqnarray*}
 &\mathbb{E}\big[f(S_{\ell_n+1}-S_{\ell_n},\ldots,S_{\ell_n+k}-S_{\ell_n}),{n-\ell_n>n\epsilon}\big]  \qquad \qquad \qquad \qquad \qquad \qquad
 \\ & 
 \geq  (\mathbb{E}_0^V[f(S_1,\ldots,S_k)]-\delta) \qquad \qquad  \qquad \qquad
 \\ & \qquad  \qquad \qquad \qquad 
\times  \sum_{i=0}^{\lceil n-n\epsilon \rceil -1} \mathbb{P}(S_j\geq S_i, \forall  j\in [ i,n])\mathbb{P}(S_j\geq S_i, \forall j\leq i) 
\\ & =  \big(\mathbb{E}_0^V[f(S_1,\ldots,S_k)]-\delta\big)\times \mathbb{P}(n-\ell_n\geq n \epsilon) .
\end{eqnarray*}
Thus, for $n>n(\delta)$,
 $$\mathbb{E}\big(f(S_{\ell_n+1}-S_{\ell_n},\ldots,S_{\ell_n+k}-S_{\ell_n})\vert n-\ell_n>n\epsilon\big)\geq \mathbb{E}_0^V[f(S_1,\ldots,S_k)]-\delta,$$
 and by the same argument,
 $$\mathbb{E}\big(f(S_{\ell_n+1}-S_{\ell_n},\ldots,S_{\ell_n+k}-S_{\ell_n})\vert n-\ell_n>n\epsilon\big)\leq \mathbb{E}_0^V[f(S_1,\ldots,S_k)]+\delta. $$
This yields the desired result.
 \qed
 \medskip

  In the above result we proved  convergence of the post-infimum process. Since the random variables $X$ are centered,
 by considering the reflected random walk $-S$, we derive a similar convergence result for the pre-infimum process to
a limit that we now introduce.
Since the environment has a density, the model enjoys a simplification: conditioning the walk to be positive is the same as 
conditioning it to be non-negative.
Define the  process $(S_k^{\downarrow})$  as the homogeneous Markov chain starting from 0  with transition function $p^{\hat{V}}$
given by \eqref{eq:Doob} with
 $$\hat{V}(x)=1+ \mathbb{E}\Big(\sum_{i=1}^{\hat{\sigma}(0)-1}1_{\{S_i\leq x\}}\Big), \quad x \geq 0,$$
 and
 $$\hat{\sigma}(0)=\min\{k\geq 1: S_k \leq 0\}.$$

 \begin{mycorol}
 \label{col:finite}
 For a bounded function $f(x_1,...,x_k)$ and $\epsilon \in (0,1)$, we have
 \begin{equation}  \nonumber
 \lim_{n\to\infty} \mathbb{E}\big(f(S_{\ell_n-1}-S_{\ell_n},
\ldots,S_{\ell_n-k}-S_{\ell_{n}})\vert \ell_n>n\epsilon\big)= \mathbb{E}(f(S_1^{\downarrow},\ldots,S_k^{\downarrow})).
 \end{equation}
 \end{mycorol}
Since the walk is centered, the event $\{n\epsilon <\ell_n<n-n\epsilon\}$ will happen with high probability for $\epsilon$ small enough. Then  Theorem \ref{theo:finite} and Corollary \ref{col:finite} imply that

 \begin{mycorol}
 \label{col:finiteCV}
  For fixed $K$, the following convergence results hold as $n\to \infty$:
   \begin{equation}  \nonumber
   (S_{\ell_n+k}-S_{\ell_n})_{1\leq k\leq K}\stackrel{\mathcal{L}}{\longrightarrow}(S_k^\uparrow)_{1\leq k\leq K}\ ,
   \end{equation}
  \begin{equation}  \nonumber
   (S_{\ell_n+k}-S_{\ell_n})_{-1\geq k\geq -K}\stackrel{\mathcal{L}}{\longrightarrow}(S_k^\downarrow)_{1\leq k\leq K}\ ,
   \end{equation}
   \begin{equation}
  \label{result:localization}
   \left(\sum_{k=-K}^K e^{-(S_k-S_{\ell_n})}\right)^{-1}\stackrel{\mathcal{L}}{\longrightarrow}
   \left(1+ \sum_{k=1}^{K} e^{-S_k^\uparrow}+\sum_{k=1}^K e^{-S_k^\downarrow}\right)^{-1}\ .
   \end{equation}
  \end{mycorol}

\subsection{Growth of random walk conditioned to stay positive}
In the literature, it is well known that the random walk conditioned to stay positive can be constructed based on an infinite number of time reversal at the ladder time set of the  walk $(S,\mathbb{P})$. The first proof is given by Golosov~\cite{Golo84}  for the case of random walk with expectation zero and later Tanaka~\cite{Tana89} gave a proof for more general case. We first present Tanaka's construction~\cite{Tana89}, and sumarize the results. 
\medskip

 Let $\{(H_k^+,\sigma_k^+)\}_{k\geq 0}$ be the sequence of strictly increasing ladder heights and times respectively of $(S,\mathbb{P})$ with $H_0^+=\sigma_0^+=0$. Define $e_1,e_2,...$ the sequence of excursions of $(S,\mathbb{P})$ from its supremum that have been time reversed:
\begin{equation}
 \nonumber
e_n=(0,S_{\sigma_n^+}-S_{\sigma_n^+-1},S_{\sigma_n^+}-S_{\sigma_n^+-2},...,S_{\sigma_n^+}-S_{\sigma_{n-1}^++1},S_{\sigma_n^+}-S_{\sigma_{n-1}^+}) \ ,
\end{equation}  
for $n\geq 1$. Write for convenience $e_n=(e_n(0),e_n(1),...,e_n(\sigma_n^+-\sigma_{n-1}^+))$ as an alternative for the step of each $e_n$. By Markov property, $e_1,e_2,...$ are independent copies of $e=(0,S_{\sigma^+}-S_{\sigma^+-1},S_{\sigma^+}-S_{\sigma^+-2},...,S_{\sigma^+}-S_{1},S_{\sigma^+})$ where $\sigma^+=\inf\{k\geq 0: S_k\in(0,+\infty)\}$ is from (\ref{eq:V}).  Tanaka's construction for the reflected random walk $(-S)$ consists in the following process $W^{\uparrow}=\{W_n^{\uparrow}: n\geq  0 \}$:
\begin{equation}
\label{params:Tanaka}
W_n^{\uparrow}= \left\{ \begin{array}{ll}
 e_1(n), & \textrm{for}\ 0\leq n\leq \sigma_1^+\\
  H_1^+ + e_2(n-\sigma_1^+) & \textrm{for}\ \sigma_1^+<n\leq \sigma_2^+ \\
  \ldots\\
  H_{k-1}^++e_k(n-\sigma_{k-1}^+) & \textrm{for} \ \sigma_{k-1}^+<n\leq \sigma_{k}^+\\
  \ldots
    \end{array} \right.
\end{equation}

Under the condition $\mathbb{P}\{\sigma^+ <\infty\}=1$, the main theorem in~\cite{Tana89} states that $\{W_n^{\uparrow}\}$ is a Markov chain process on $[0,+\infty)$ with transition function $\hat{p}(x,dy)$, which is given by 
 \begin{equation}
  \nonumber
 \hat{p}(x,dy)= \frac{g(y)}{g(x)} \mathbb{P}(x+X\in dy) 1_{(0,+\infty)}(y) \ ,
 \end{equation}
 where $$g(x)=\mathbb{E}\Big[\sum_{n=0}^{\sigma^+}1_{\{-x < S_n\}}\Big].$$
 As we consider here  log-gamma variables with $\theta=\mu/2$, then we have $\mathbb{P}$-a.s $\sigma^+=\sigma(0)<1$ and $g=V$ from (\ref{eq:V}). Therefore the process $W^{\uparrow}$ has the same law as the random walk conditioned to stay non negative, i.e., $S^{\uparrow}$ defined above proposition \ref{theo:convergence2 }.
This identity provides an elegant way to determine  the growth rate of the limit process $S^{\uparrow}$.

  \begin{mylem} \label{lem:growth0}
  For every $\epsilon >0$, then we have :
  \begin{equation}
   \nonumber
  \lim_{n\to\infty} \frac{S_n^\uparrow}{n^{1/2-\epsilon}} = +\infty, \ \mathbb{P} -a.s  \ .
  \end{equation}
  As a consequence, for fixed $\delta > 0$ there exists a constant $k=k(\delta)$ such that 
  \begin{equation}  \nonumber
  \mathbb{P}(S_n^{\uparrow }\geq n^{1/2- \epsilon},\forall \ n\geq k)\geq 1-\delta \ .
  \end{equation}
  
  \end{mylem}
 \textit{ Proof:} We follow the lines of \cite{HaKeKy03}.
 Let $\{M_k^+,v_k^+\}_{k\geq 0}$ be the space-time points of increase of the future minimum of $(W^{\uparrow},\mathbb{P})$. That is, $M_0^+=v_0^+=0$,
 \begin{equation}
  \nonumber
 v_k^+= \inf\{n>v_{k-1}^+: \min_{r\geq n}W_r^{\uparrow}=W_n^{\uparrow}\} \ \textrm{and} \ M_k^+ = W_{v_k^+}^{\uparrow},
 \end{equation}
 for $k\geq 1$. From the construction of $W^{\uparrow}$, we can deduce that for each path, the sequence $\{(M_k^+,v_k^+)\}_{k\geq 0}$ corresponds precisely to   $\{(H_k^+,\sigma_k^+)\}_{k\geq 0}$. Let $L=\{L_n\}_{n\geq 0}$ be the local time at the maximum in $(S,\mathbb{P})$, that is 
 \begin{equation}
  \nonumber
 L_n=|\{k\leq n : \max_{i\leq k}S_i= S_k\}|  \ .
 \end{equation}
 Because $W^{\uparrow}$ is obtained by time reversal from $S$, then $L$ is also the local time at the future minimum of $(W^{\uparrow},\mathbb{P})$. Hence it's clear that :
 \begin{equation}
  \nonumber
 S_n\leq H_{L_n}^+=M_{L_n}^+\leq W_n^{\uparrow}.
 \end{equation} 
 Now we need the following lemma (e.g., Theorem 3 in~\cite{HaKeKy03}):
 \begin{mylem}\label{lem:growth}
 Consider the random walk $(S,\mathbb{P})$. Now suppose that $\Phi \downarrow 0$ and that $\mathbb{E}(S_1)=0$ and $\mathbb{E}(S_1^2)<\infty$. Then 
 \begin{equation}  \nonumber
 \mathbb{P}_x(\max_{k\leq n}S_k < \sqrt{n}\Phi(n)\ i.o.)=0\ or\ 1 \ .
 \end{equation}
 according to
 \begin{equation*}
 \int_1^{\infty} \frac{\Phi(t)}{t} dt < \infty \ or = \infty  \ .
 \end{equation*}
 \end{mylem}
 We use the standard notations  "i.o." for "infinitely often" and   "ev." for "eventually". 
 For $\Phi(n) =n^{-\epsilon} $,  the integral converges and Lemma  \ref{lem:growth} yields
 \begin{eqnarray}  \nonumber
   1 =& \mathbb{P}(\max_{i\leq n}S_i \geq \sqrt{n}\Phi(n) \ \text{ev}. )\\  \nonumber
 =& \mathbb{P}(W_n^{\uparrow}\geq  \sqrt{n}\Phi(n) \ \rm{ev}. )   \ .\end{eqnarray}
Again using the fact that we may replace $\Phi$ by $c\Phi$ for any $c>0$, and that the integral in the lemma still converges, it follows easily that
 \begin{equation}
 \nonumber 
 \liminf_{n\to\infty} \frac{W_n^{\uparrow}}{ \sqrt{n}\Phi(n)} = \infty \ ,
 \end{equation}
 $\mathbb{P}$-almost surely. As $W^{\uparrow}$ and $S^{\uparrow}$ have the same law under $\mathbb{P}$, we get the first result. Then it's clear for fixed $\delta$, there exists $k$ such that 

 \begin{equation*}
\mathbb{P}(S_n^{\uparrow}\geq n^{1/2- \epsilon},\forall \ n\geq k)=\mathbb{P}(A_k)\geq 1-\delta  \ .
\end{equation*}
\qed
\medskip

We complement Lemma \ref{lem:growth0} with the following 
 version for the conditioned random walk, proved by Ritter~\cite{Ritter81}. 
 \begin{mytheo}[\cite{Ritter81}]\label{th:ritter}
 Fixed $\eta< 1/2$ then :
 \begin{equation}
  \nonumber
 \lim_{\delta\to 0}  \inf_n\mathbb{P}[\inf_{k\leq n}(S_k-\delta k^\eta)\geq 0|\tau>n]=1\ .
 \end{equation}
 \end{mytheo}
Note now that the time $\ell_n$ of the first  minimum of $S$ on $[0,n]$  is such that,  for fixed $\epsilon \in (0,1)$ we have for large $n$,
 \begin{equation}  \label{lemma:ell_n}
 \mathbb{P}[\ell_n<(1-\epsilon^2)n]>1-\epsilon\ .
 \end{equation}
Indeed, by  the invariance principle of Donsker(1951), $\ell_n/n$ converges in law to the time of the global minimum of the standard Brownian motion on $[0,1]$, which obeys the arcsine distribution~\cite[problem 8.18]{KaSh98}.
Therefore, conditionally on  the event $\{\ell_n<(1-\epsilon^2)n\}$,   Theorem \ref{th:ritter} gives us the growth of the random walk after  the minimum:
 \begin{mycorol}\label{Cor:Ritter}
 If $\eta \in (0,1/2)$, then uniformly in $n$ :
 \begin{equation}  \nonumber
 \lim_{\delta\to 0} \mathbb{P}[S_{k+\ell_n}-S_{\ell_n}>\delta k^\eta\  \textrm{for all} \ k\leq n-\ell_n] = 1\ .
 \end{equation}
 \end{mycorol}

 \textit{ Proof:}
 The proof is similar to the proof of Proposition \ref{theo:finite}. To simplify the notation define
 \begin{equation*}
 A_{\delta}=\{S_{k+\ell_n}-S_{\ell_n}>\delta k^\eta\  \textrm{for all} \ k\leq n-\ell_n\}\ ,
 \end{equation*}
 and
 \begin{equation*}
 A_{\delta,j}=\{S_{k+\ell_n}-S_{\ell_n}>\delta k^\eta\  \textrm{for all} \ k\leq n-\ell_n,\ \ell_n=j\}\ .
 \end{equation*}
 Then we have
 \begin{equation*}
\mathbb{P}[A_{\delta}]=\sum_{j=1}^{n}\mathbb{P}[A_{\delta,j}]=
\sum_{j=1}^{n}\mathbb{P}[A_{\delta,j}\vert \ell_n=j]\mathbb{P}[\ell_n=j]\ .
 \end{equation*}
 We know that the event  $\{\ell_n=i\}$ can be written as
 $$\{\ell_n=i\}=\{S_j\geq S_i, \forall j\leq i\}\cap \{S_j\geq S_i, \forall  j\in [ i,n]\}\ ,$$
 Both random variables $ A_{\delta,j}$ and $1_{\{S_j\geq S_i, \forall  j\in [ i,n]\}}$ are measurable with respect to $\sigma(X_{i+1},...,X_n)$ and are independent of  the event $\{S_j\geq S_i, \forall j\leq i\}$. By the Markov property, it follows that
 \begin{equation*}
 \mathbb{P}[A_{\delta,j}\vert \ell_n=j]=\mathbb{P}[S_{k}>\delta k^\eta\  \textrm{for all} \ k\leq n-j\vert \tau>n-j]\ .
 \end{equation*}
 For $\epsilon>0$ by using Theorem \ref{th:ritter}, there exists $\delta_{\epsilon}$ and $n_{\epsilon}$ such that for for all $\delta<\delta_{\epsilon}$, $m> n_{\epsilon}$ 
 \begin{equation*}
\mathbb{P}[\inf_{k\leq m}(S_k-\delta k^\eta)\geq 0|\tau>m]\geq 1-\epsilon\ .
 \end{equation*}
 By putting $m=n-j$ and summing $j\in \{1,2,...,[(1-\epsilon )n]\}$ we obtain 
 \begin{equation*}
\mathbb{P}[A_{\delta}]\geq (1-\epsilon)\mathbb{P}[\ell_n<(1-\epsilon)n]\ .
 \end{equation*}
 So by (\ref{lemma:ell_n}), for $n$ large enough, we have
 \begin{equation*}
\mathbb{P}[A_{\delta}]\geq 1-2\epsilon\ ,
 \end{equation*}
 which implies easily the corollary.
\qed

\section{Proof of the main results in the point-to-line case} \label{section: proof}
We split the section according to $\theta=\mu/2$ or not, starting with the first case, which is more involved than the second one.
The reason why $\theta=\mu/2$ is special is that the random variable $X$ in (\ref{eq:RW}) is 
centered, and even symmetric.

\subsection{Equilibrium case}\label{section: proof-eq}
In the equilibrium setting $\theta=\mu/2$, we know that the post- and pre-infimum chain converge in law to the random walk conditioned to stay positive. As these limit processes grow fast enough, we can indeed prove that the endpoint densities of the polymer converge when its length goes to infinity. Firstly we  consider the distribution at the favourite endpoint, and we later extend the arguments to all the points: 
\begin{mylem}\label{lem:one-point} 
For $n\to\infty$,
\begin{equation}
\nn
\xi_{\ell_n}^n=\left( \sum_{i=0}^n e^{-(S_i-S_{\ell_n})}\right)^{-1} \stackrel{\mathcal{L}}{\longrightarrow}\xi_0 =\left( 1+\sum_{i=1}^\infty e^{-S_i^\uparrow}+\sum_{i=1}^\infty e^{-S_i^\downarrow}\right)^{-1} \ .
\end{equation}
\end{mylem}

\textit{Proof.}
From Lemma \ref{lem:growth0}, the random walk conditioned to stay positive is lower bounded by some factor of $n^{1/2-\epsilon}$, thus the random variable $\xi_0$ is well defined and strictly positive.
By the continuous mapping theorem, the claim is equivalent to convergence in law of the inverse random variables.
Then, in order to prove the lemma,
 it suffices to show that, for all bounded and uniformly continuous function $f$, we have
\begin{equation}
\label{proof:convegence}
\mathbb{E}\Big[f\Big(\sum_{i=0}^n e^{-(S_i-S_{\ell_n})}\Big )\Big]\longrightarrow \mathbb{E}\Big [f\Big(1+\sum_{i=1}^\infty e^{-S_i^\uparrow}+\sum_{i=1}^\infty e^{-S_i^\downarrow}\Big)\Big]\  ,
\end{equation}
as $n \to \8$. By  (\ref{result:localization}) in Corollary \ref{col:finiteCV}, we already know that,  for a fixed $K$,
\begin{equation}
\label{proof:convergence1}
\mathbb{E}\Big[f\Big(\sum_{i=\ell_n-K}^{i=\ell_n+K} e^{-(S_i-S_{\ell_n})}\Big )\Big]\longrightarrow \mathbb{E}\Big[f\Big(1+\sum_{i=1}^K e^{-S_i^\uparrow}+\sum_{i=1}^K e^{-S_i^\downarrow}\Big)\Big]\ .
\end{equation}
By uniform continuity, 
given an $\epsilon >0$ there exists $\rho>0$ such that $|f(x)-f(y)|<\epsilon$ for all $x,y$ with $|x-y|< \rho$. Now we will prove that, for all positive $\epsilon$, we can find finite $K=K(\epsilon)$ and $n_0(\epsilon)$ such that for $n \geq n_0(\epsilon)$, 
it holds
\begin{equation}
\label{proof:convergence2}
\Big| \mathbb{E}\Big[f\Big(\sum_{i=0}^n e^{-(S_i-S_{\ell_n})}\Big )\Big]-\mathbb{E}\Big[f\Big(\sum_{i=\ell_n-K}^{i=\ell_n+K} e^{-(S_i-S_{\ell_n})}\Big )\Big]\Big|<\epsilon\ ,
\end{equation}
and
\begin{equation}
\label{proof:CV4}
\Big \vert \mathbb{E}\Big[f\Big(1+\sum_{i=1}^K e^{-S_i^\uparrow}+\sum_{i=1}^K e^{-S_i^\downarrow}\Big)\Big]-\mathbb{E}\Big [f\Big(1+\sum_{i=1}^\infty e^{-S_i^\uparrow}+\sum_{i=1}^\infty e^{-S_i^\downarrow}\Big)\Big] \Big \vert<\epsilon\ .
\end{equation}
Then, by combining  $(\ref{proof:convergence1})$, $(\ref{proof:convergence2})$ and $(\ref{proof:CV4})$ we get $(\ref{proof:convegence})$ and the proof is finished.

In order to get  $(\ref{proof:convergence2})$, it is enough to prove that for all positive $\rho, \epsilon$ there exists a finite $K$ such that
\begin{equation}
\label{proof:convergence2b}
\mathbb{P}\Big ( \sum_{i=0}^{\ell_n-K}e^{-(S_i-S_{\ell_n})}+\sum_{i=\ell_n+K}^{n}e^{-(S_i-S_{\ell_n})}
< \rho \Big) > 1-\epsilon
\end{equation}
for all $n$ large enough, while,   
 in order to get  $(\ref{proof:CV4})$, it is enough to prove that for all positive $\rho, \epsilon$ there exists a finite $K$ such that
\begin{equation}
\label{proof:CV4b}
\mathbb{P}\Big ( \sum_{i \geq K}e^{-S_i^\uparrow}+  \sum_{i \geq K}e^{-S_i^\downarrow}
< \rho \Big) > 1-\epsilon
\end{equation}
for all $n$ large enough.

To prove $(\ref{proof:convergence2b})$, we use Corollary \ref{Cor:Ritter}: For any fixed $\eta<1/2$, 
 choose $\delta>0$, such that for all $n\in\mathbb{N}$ :
\begin{equation*}
\mathbb{P}[S_{k+\ell_n}-S_{\ell_n}>\delta k^\eta;  k=1,2,\ldots, n-\ell_n] \geq  1-\epsilon/2\ .
\end{equation*}
Because the random variable $X$ is symmetric,
 the pre-infimum process verifies the same properties, i.e for $n\in \mathbb{N}$ :
\begin{equation*}
\mathbb{P}\big[S_{\ell_n-k}-S_{\ell_n}>\delta k^\eta; k=1,\ldots, \ell_n\big] \geq  1-\epsilon/2\ .
\end{equation*}
Then, choosing $K$ such that 
$\sum_{k=K}^\infty e^{-\delta k^\eta}< \rho/2$ yields $(\ref{proof:convergence2b})$.
A similar argument  leads to (\ref{proof:CV4b}).
This completes the proof of the lemma.
\qed
\medskip

In the course of the proof we have discovered the limit endpoint densities $(\xi_k)_{k\in Z}$ is given by formula \eqref{eq:thelimit}.
Repeating the argument in the proof of Lemma \ref{lem:one-point}, 
it is straightforward to extend the result  to a finite set of points around the maximum point:  
\begin{mylem}\label{lem:finite-points}
For fixed K, $n\to\infty$ :
\begin{equation*}
(\xi_{\ell_n+k}^n)_{-K\leq k\leq K} \stackrel{\mathcal{L}}{\longrightarrow}(\xi_k)_{-K\leq k\leq K}
\end{equation*}
\end{mylem}
\textit{Proof of Theorem \ref{theo:main} in the case of  $\theta=\mu/2$.} 
Recall that $\tilde \xi^{(n)}$ and $\xi^n$ have the same law, and that the total variation distance between probability measures on $\mathbb Z$ coincides with the 
$\ell^1$-norm.
Taking a function $f:\ell^1\rightarrow R$ bounded and uniformly continuous in the norm $\vert.\vert_1$ and one needs to prove that
\begin{equation}\label{proof:ar1}
\mathbb{E}[f((\xi_{\ell_n+k}^n)_{k\in \mathbb{Z}})] \rightarrow \mathbb{E}[f((\xi_k)_{k\in \mathbb{Z}})].
\end{equation}
We will use almost the same idea as in the proof of Lemma \ref{lem:one-point}. 
From Lemma \ref{lem:finite-points} we have, with a slight abuse of notation,
 \begin{equation}\label{proof:ar22} 
\mathbb{E}[f((\xi_{\ell_n+k}^n)_{k\in [-K,K]})] \rightarrow \mathbb{E}[f((\xi_k)_{k\in  [-K,K]})]
\end{equation}
Fixing $\epsilon>0$, there exists by continuity  some $\delta>0$ such that $|x-y|_1<\delta$ implies $|f(x)-f(y)|<\epsilon$. Hence,
\begin{equation}\label{proof:ar3}
\mathbb{E}\Big[f((\xi_{\ell_n+k}^n)_{k\in [-K,K]})-f((\xi_{\ell_n+k}^n)_{k\in \mathbb{Z}})\Big]\leq\epsilon
\end{equation}
provided that
\begin{equation} \label{proof:remark}
\mathbb{E}[\sum_{k:|k|>K}\xi_{\ell_n+k}^n ]\leq \delta\ ,
\end{equation}
and similarly for $\xi$ instead of $\xi^n$.
Since $\mathbb{E}[\xi_k]$ is a probability measure on $\mathbb{Z}$, we can take $K$ large enough so that $\mathbb{E}[\sum_{k:|k|>K}\xi_k ]\leq \delta$. 
Then, from Lemma \ref{lem:finite-points}, we see that, as $n \to \8$, 
$$\mathbb{E} \big[\sum_{k:|k|>K}\xi_{\ell_n+k}^n \big]=1-\mathbb{E}\big[\sum_{k:|k|\leq K}\xi_{\ell_n+k}^n \big] \longrightarrow \mathbb{E}\big[\sum_{k:|k|>K}\xi_k \big]
=1-\mathbb{E}\big[\sum_{k:|k|\leq K}\xi_{k} \big] ,$$ yielding (\ref{proof:remark}).  By combining (\ref{proof:ar22}) and (\ref{proof:ar3}), we obtain (\ref{proof:ar1}).
\qed 
\medskip

\textit{Proof of Corollary \ref{theo:loc}.}
It is enough to note that
$$
 \lim_{n\to \infty} Q_n^{\omega}[\vert x_n \cdot {\bf e_1}-l_n\vert \geq K] = \sum_{k:|k|>K}\xi_k \quad \textmd{ in law},
$$
which vanishes as $K \to \8$. \qed
\medskip

\textit{Proof of Corollary \ref{cor:first}.}
Recall first that under $Q_{n-1}^\omega$, the steps after the final time $n-1$ are uniformly distributed, and  independent from everything else.
Then, it is enough to note that
\begin{eqnarray*}
Q_{n-1}^{\omega}[x_n=x] &=& 
Q_{n-1}^{\omega}[x_{n-1}=x-{\mathbf e_1}, x_{n}-x_{n-1}={\mathbf e_1}]  + \\
&& \qquad \qquad
 Q_{n-1}^{\omega}[x_{n-1}=x-{\mathbf e_2}, x_{n}-x_{n-1}={\mathbf e_2}] \\
 &=& \Big( Q_{n-1}^{\omega}[x_{n-1}=x-{\mathbf e_1}]  + 
 Q_{n-1}^{\omega}[x_{n-1}=x-{\mathbf e_2}] \Big)\times
 \frac{1}{2} \\
 & \stackrel{\cal L}{\longrightarrow}&  \frac{\xi_k+\xi_{k+1}}{2},
\end{eqnarray*}
with $k$ determined by $(x-{\mathbf e_1})\cdot {\mathbf e_1}= l_{n-1}+k$. \qed
\medskip

\noindent
Now we give the proof for Theorem \ref{theo:arcsin} and \ref{theo:deviation}:
\medskip

\textit{Proof of Theorem \ref{theo:arcsin} for $\theta=\mu/2$.} First, recall that $l_n$ and $\ell_n$ have the same law, so we can focus on the latter one. 
By definition of $\ell_n$ and  Donsker's invariance principle, we have directly 
\begin{equation*}
\frac{\ell_n}{n}=\frac{1}{n} \underset{i\in[0,n]}{\arg\min}\  S_i \stackrel{\mathcal{L}}{\longrightarrow}\underset{t\in[0,1]}{\arg\min}\ W_t\ .
\end{equation*}

By  L\'evy's arcsine law~\cite{Levy39}, the location of the minimum  of the Brownian motion,
i.e. the above limit, 
has the density $\pi^{-1}\big(s(1-s)\big)^{-1/2}$
on the interval $[0,1]$.\qed
\medskip

\textit{Proof of Theorem \ref{theo:deviation} for $\theta=\mu/2$.}
We can express the first term in (\ref{equa:deviation}) as 
\begin{equation} 
\label{equa:decom}
\frac{-1}{\sqrt{n}}\log Q_n^{\omega}\big(x_n=([ns],n\!-\![ns])\big)
\stackrel{\cal L}{=} \frac{-1}{\sqrt{n}}(S_{[ns]}-S_{\ell_n})+\frac{1}{\sqrt{n}}\log(\sum_{k=1}^{n}e^{-(S_k-S_{\ell_n})})\ ,
\end{equation}
with $\stackrel{\cal L}{=} $ the equality in law.
As the second term in the right-hand side
 is almost surely dominated by $\frac{\log n}{\sqrt{n}}$, then again,  Donsker's invariance principle yields (\ref{equa:deviation}).

On the other hand, let A be an interval.  For all $s\in A$, we have 
\begin{eqnarray*}
Q_n^{\omega}\big(x_n=([ns],n\!-\![ns])\big) &\leq &Q_n^{\omega}(x_n\in nA) \\
&\leq& n\max_{x\in A}Q_n^{\omega}\big(x_n=([nx],n\!-\![nx])\big)\ ,
\end{eqnarray*}
which  means that
\begin{eqnarray*}
\max_{x\in A}Q_n^{\omega}\big(x_n=([ns],n\!-\![ns])\big)
&\leq & Q_n^{\omega}(x_n\in nA)\\
&\leq & n\max_{x\in A}Q_n^{\omega}\big(x_n=([nx],n\!-\![nx])\big)\ .
\end{eqnarray*}
From Donsker's invariance principle it follows 
\begin{eqnarray*}
\frac{-1}{\sqrt{n}}\max_{x\in A} \log Q_n^{\omega}\big(x_n\!=\!([ns],n\!-\![ns])\big)&\stackrel{\cal L}{=} & \min_{x\in A} \frac{S_{[nx]}\!-\!S_{\ell_n}}{\sqrt{n}} + \frac{1}{\sqrt{n}}\log(\sum_{k=1}^{n}e^{-(S_k\!-\!S_{\ell_n})}) \\
&\stackrel{\mathcal{L}}{\longrightarrow}& \min_{x\in \overline{A}}W(x) -\min_{[0,1]} W\ ,
\end{eqnarray*}
which, in turn,  yields (\ref{equa:interval}). \qed

\subsection{Non-equilibrium case} \label{sec:4.2}

\textit{Proof of Theorem \ref{theo:main} in the case of  $\theta\neq\mu/2$.}
Without loss of generality, we assume that $\theta< \mu/2$, which implies that $m=\mathbb{E}[X]>0$ and the random walk $S$ drifts to $+\infty$. By the law of large number, we have for all $a \in (0,m)$,
$$
M=\min_n (S_n-na)>-\infty \quad \mathbb{P}-a.s
$$
It follows that $\mathbb{P}$-a.s for every integer $n$,
$$
e^{-S_n}\leq e^{-na-M}
$$
Then the sum of $e^{-S_n}$ converges $\mathbb{P}$-a.s and we can identify the limit distribution $\xi$ as
$$
\overline{\xi}_k=\frac{e^{-S_k}}{\sum_{i=0}^{\infty}e^{-S_i}}\,.
$$
Indeed, it is clear that, for $k\in \mathbb{Z}_+$,
$$
\xi_k^n= \frac{e^{-S_k}}{\sum_{i=0}^{n}e^{-S_i}}\to \overline{\xi}_k ,\ \ \ \mathbb{P}-a.s
$$
Since the random walk drifts to $+\infty$   the global minimizer 
$$
\ell=\arg\min\{S_k,k\in\mathbb{Z}_+\}
$$
is $\mathbb{P}$-a.s  finite. 
Moreover, for $n$ large enough we have $\ell_n=\ell$,  and by centering the measure $\xi^n$ and $\overline{\xi}$ around $\ell_n$ and $\ell$ respectively, we can easily obtain that

\begin{equation} \nonumber
 \tilde \xi^{(n)} \stackrel{\mathcal L}{\longrightarrow} \xi \qquad \textmd{ in the space } ({\mathcal M}_1, \| \cdot \|_{TV}),
\end{equation}
where  $\tilde \xi^{(n)}$ is defined as in (\ref{eq:cvloixi}) and $\xi_k=\overline{\xi}_{\ell+k}$. This yields Theorem \ref{theo:main} in the case $\theta<\mu/2$.
\qed
\medskip 

\textit{Proof of Theorem \ref{theo:arcsin} for $\theta \neq \mu/2$.} It is a straightforward consequence of the above, since $\ell_n={\mathcal O}(1)$ if $\theta < \mu/2$ or 
 $n-\ell_n={\mathcal O}(1)$ if $\theta > \mu/2$. \qed
 
\medskip 

 \textit{Proof of Theorem \ref{theo:deviation} for $\theta\neq \mu/2$.} Though it was already proved in~\cite{Sepp12, GeorSepp13}, we give another argument for completeness. Applying the law of large numbers for i.i.d. variables in (\ref{equa:decom}), we directly obtain the claim. \qed


\section{ Localization of the point-to-point measure} \label{sec:last}
 In this section, we consider the point-to-point measure with mirror boundary conditions. Recall the definition 
 of the model P2P-b.c.$(\Theta)$ from (\ref{eq:P2Pbc}).
  
 In this situation, beside the usual partition function $Z_{m,n}$,  we will also define the reverse partition function $\tilde{Z}_{m,n}$ for $(m,n) \in R_N$ as
 \begin{equation} \label{eq:1aa}
 \tilde{Z}_{m,n}= \sum_{\bx \in \tilde{\Pi}_{m,n}^N }\prod_{t=m+n }^{(p+q)N-1} Y_{x_t},
 \end{equation}
 where $\tilde{\Pi}_{m,n}^N$ denotes the collection of up-right paths $\bx=(x_t; m+n\leq t \leq (p+q)N)$ in the rectangle $R_N$ that go from $(m,n)$ to $(pN,qN)$.  Note that in the reverse partition function we exclude the weight at $(pN,qN)$. 
Also, it depends on $N, p$ and $q$, but we  omit to indicate it in the notation. 
Moreover, we can also define the ratios $\tilde{U}$ and $\tilde{V}$ as in the usual case,
 \begin{equation} \label{eq:1a}
 \tilde{U}_{m,n}= \frac{\tilde{Z}_{m,n}}{\tilde{Z}_{m+1,n}}
 \end{equation}
\begin{equation} \label{eq:1b}
 \tilde{V}_{m,n}= \frac{\tilde{Z}_{m,n}}{\tilde{Z}_{m,n+1}}
 \end{equation}
 If we take the point $(pN,qN)$ as the initial point then the reverse environment $(\tilde{Z},\tilde{U},\tilde{V})$ is also a stationary log-gamma system with boundary conditions.  Indeed one sees from (\ref{eq:1aa}), (\ref{eq:1a}) and (\ref{eq:1b}) that $ \tilde{U}_{m,qN}=Y_{m,qN}$ and $ \tilde{V}_{pN,v}=Y_{pN,n}$.
\begin{figure}[htbp]\label{des1}
\hspace*{-2cm}
\centering
\includegraphics[scale=0.7,width=13cm,height=5.3cm]{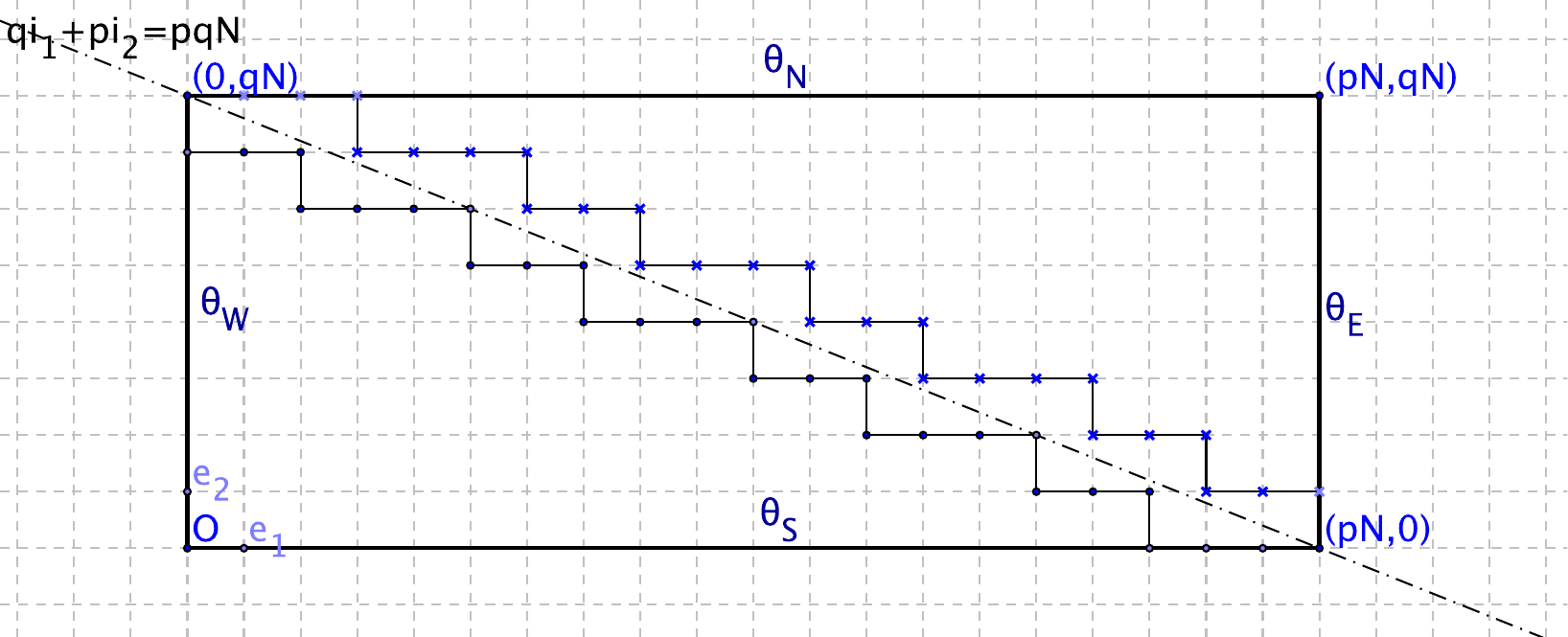}  
\caption{ Upper and lower "transverse diagonal" with $p=5, q=2, N=4$. Their vertices are indicated by dots and crosses respectively. $H_{p,q}^{-,N}$ is the region below the diagonal. The boundary conditions are indicated on the boundaries of the rectangle.}
\end{figure}
 We partition the rectangle $R_N$ according to the lower half space  
 $$H_{p,q}^{-,N}=\{(i_1,i_2)\in \mathbb{Z}^2: qi_1+pi_2\leq pqN\}$$
and its complement. 
 		   In order to simplify the notations, we denote for an edge $f$ with endpoints in  $R_N$, 
 		   \begin{equation} \label{eq:1c}
 		    T_f= \left\{ \begin{array}{ll}
 		    U_x & \textrm{if} \ f\in H_{p,q}^{-,N} \ \textmd{and}\  f =\langle \cdot ,x\rangle \ \textmd{is horizontal}\\
 		      V_x^{-1} & \textrm{if} \ f\in H_{p,q}^{-,N}\ \textmd{and}\  f =\langle  \cdot ,x\rangle \ \textmd{is vertical}  \\
 		      \tilde{U}_x^{-1}  & \textrm{if} \ f\notin H_{p,q}^{-,N}\ \textmd{and}\  f =\langle x, \cdot \rangle \ \textmd{is horizontal}  \\
 		       \tilde{V}_x  & \textrm{if} \ f\notin H_{p,q}^{-,N},\ \textmd{and}\  f =\langle x, \cdot \rangle \ \textmd{is vertical} \\
 		       \end{array} \right.
 		   \end{equation}
 From 
Fact 1 of section~\ref{subsec:p2l} 
and independence of the weights in $H_{p,q}^{-,N}$ and 
 its complement, for every down-right path $z$ it follows that the variables $\{T_f:f\in z\}$ are mutually independent.
 We stress that independence relies also on the expressions of $Z_{m,n}$ and  $\tilde{Z}_{m',n'}$, where there are no shared weights. 
 The marginal distribution of $T_f$ is given by the stationary structure, it is a log-gamma distribution with the  appropriate parameter.
 Let $\partial H_{p,q}^{-,N}$ be the transverse diagonal in $R_N$, which is given as 
 	  $$\partial H_{p,q}^{-,N}=\{(i_1,i_2)\in H_{p,q}^{-,N}: (i_1+1,i_2+1)\notin H_{p,q}^{-,N}\}.$$	
 Consider  the "lower transverse diagonal" given as
  $$\mathcal{L}_{p,q}^N  \  \textmd{is the down-right  path} \ \bx=(x_i): (0,qN)\to(pN,0) \ \textmd{ with } x_i\in \partial H_{p,q}^{-,N} ,$$
 and the "upper transverse diagonal",
 $$\mathcal{U}_{p,q}^N=(1,1)+\mathcal{L}_{p,q}^N,$$
 see Figure 
 2. 
 Define also the set of up-right edges across the transverse diagonal,
  $$\mathcal{A}_{p,q}^N  =\{\langle z_1,z_2\rangle: z_1\in H_{p,q}^{-,N}, z_2 \notin H_{p,q}^{-,N}, \vert z_1- z_2\vert_1=1\}.$$
 Each up-right path $\bx$ that goes from $(0,0)$ to $(pN,qN)$,  intersects the transverse diagonal once and only once. Precisely, the mapping 
\begin{equation} \nonumber
 \bx=(x_j)_j  \mapsto G(\bx)=\langle x_i, x_{i+1}\rangle , \textmd{ with }  \langle x_i, x_{i+1}\rangle \in \mathcal{A}_{p,q}^N , 
\end{equation}
is well defined, and it indicates where the crossing takes place. (We have $i=t^-$ in (\ref{eq:transv}).)
Our main question in this section is the behaviour of the crossing edge when $N$ increases. 
By definition of the polymer measure,  for $\langle z_1,z_2\rangle\in \mathcal{A}_{p,q}^N$ we can write, with the notation $\sum_*$ for 
the sum over $\bx \in \Pi_{pN,qN}, x_{t^-}=z_1, x_{t^-+1}=z_2$,
\begin{eqnarray*}
Q_{pN,qN}^{\omega}(G(\bx)=\langle z_1,z_2\rangle) &=&   \frac{1}{Z_{pN,qN}} \;
\sum_*
\exp\Big\{\sum_{t=1}^{t^-}\omega(x_t)+\sum_{t=t_-+1}^{(p+q)N}\omega(x_t)\Big\} \\
&=&   \frac{1}{Z_{pN,qN}} \; Z_{z_1} \times \tilde{Z}_{z_2} \times \exp \{ \omega(pN,qN)\}\;,
\end{eqnarray*}
where the last factor is the contribution of the last point $x_{(p+q)N}=(pN,qN)$. In view of  (\ref{eq:1a}), (\ref{eq:1b}), (\ref{eq:1c}), this term can be expressed as
\begin{eqnarray*}
Q_{pN,qN}^{\omega}(G(\bx) =\langle z_1,z_2\rangle) 
&=&\frac{1}{Z_{pN,qN}} \; Z_{0,qN} \times \tilde{Z}_{1,qN} \times \\
&& \qquad \exp\big(-\sum_{\Pi_1(z_1) }\log(T_f) -\sum_{\Pi_2(z_2)} \log(T_f)\big) 
\\ &&  \qquad  \qquad  \times
\exp \{ \omega(pN,qN)\}\;,
\end{eqnarray*}
 where $\Pi_1(z_1)$ is the restriction of $ \mathcal{L}_{p,q}^N $ from $(0,qN)$ to $z_1$, 
and $\Pi_2(z_2)$  is the restriction of $\mathcal{U}_{p,q}^N $ from $(1,qN)$ to $z_2$. 
Note that, when computing the ratio of the left-hand side for two different values of the crossing edge $\langle z_1,z_2\rangle$, both the first and last lines 
of the right-hand side cancel.
 Thus, we consider
 $$W(\langle z_1,z_2\rangle)=\sum_{\Pi_1(z_1)}\log(T_f) +\sum_{\Pi_2(z_2)} \log(T_f)\;.$$ 
Observe that the variables $\{T_f: f\in \Pi_1(z_1)\}$
and  $\{T_f: f\in \Pi_2(z_2)\}$ are independent but not  identically distributed, so that in order to apply the same method as in previous section, we should divide $\mathcal{U}_{p,q}^N$ and $\mathcal{L}_{p,q}^N$ into identical blocks to obtain a centered random walk. Blocks are shifts of the right triangle with vertices ${\bf 0},
p{\bf e_1}, q{\bf e_2}$.
  Precisely  we denote by
   $$z_1^k=(kp,(N-k)q), \  \textmd{for}\  0\leq k\leq N, $$
the vertices in $R_N$ which sit on the line of equation $qi_1+pi_2=pqN$,  by
    $$z_2^k=z_1^k+{\bf e_1},$$
    and by $\mathcal{A}$  the set of crossing edges in the basic block $R_1$ shifted by $q{\bf e_2}$,
    $$\mathcal{A}=\Big\{\langle z_1,z_2\rangle \textmd{ up or right edge:} \ z_1+ q{\bf e_2}
    \in \mathcal{L}_{p,q}^1, \ z_2+q{\bf e_2} \in \mathcal{U}_{p,q}^1 , 1 \leq  z_2 \cdot {\bf e_1} \leq p\Big\}$$
  as shown in Figure 3.
   Note that $\langle {\bf 0},{\bf e_1}\rangle\in \mathcal{A}$ and the shifted edge $(p,-q)+ \langle {\bf 0},{\bf e_2}\rangle\in \mathcal{A}$ but $\langle {\bf 0} ,{\bf e_2}\rangle\notin \mathcal{A}$. We will use $\mathcal{A}$ and the set $(\langle z_1^k,z_2^k\rangle)_k$ to parametrize the set $\mathcal{A}_{p,q}^N$ as follows.
   For $\langle z_1,z_2\rangle \in \mathcal{A}_{p,q}^N$ we can find a unique $k$ such that, relative to any coordinate, 
   $z_1$ is between $z_1^k$ and $z_1^{k+1}$.  Then by translation, there exists a unique edge $a\in \mathcal{A}$ such that:
   $$\langle z_1,z_2\rangle =\langle z_1^k,z_2^k\rangle +a$$
   
   \begin{figure}[htbp]
\hspace*{-0cm}
\centering
\includegraphics[scale=0.7,width=10cm,height=5cm]{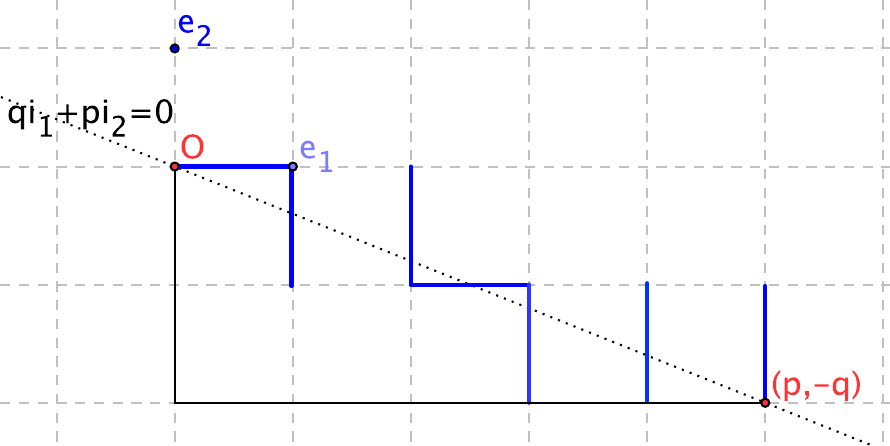}        
\caption{ The set $\mathcal A$  with $p=5, q=2$ contains 7 crossing edges in the first block,
indicated by solid lines.}
\label{Fig:2}
\end{figure}
   
   On the other hand, we have:
    
  \begin{equation} 
   \nonumber
 W_k=W(\langle z_1^k,z_2^k\rangle )=\sum_{\Pi_1(z_1^k)}\log(T_f) +\sum_{\Pi_2(z_2^k)} \log(T_f)=\sum_{i=1}^{k} X_i
\end{equation}
  where 
    \begin{equation}  \label{eq:X}
X_k= \sum_{f\in \Pi_1(z_1^k,z_1^{k+1})} \log(T_f)+\sum_{f\in \Pi_2(z_2^k,z_2^{k+1})} \log(T_f)
\end{equation}
 with $\Pi_1(z_1^k,z_1^{k+1})$ the restriction of $\mathcal{L}_{p,q}^N$  from $z_1^k$ to $z_1^{k+1}$ and  $\Pi_2(z_2^k,z_2^{k+1})$ the restriction  of $\mathcal{U}_{p,q}^N$  from $z_2^k$ to $z_2^{k+1}$.

Let $\mathcal{B}=\{ \textmd{edges} \ f \in \Pi_1(z_1^1)\cup\Pi_2(z_2^1)\}$, the edge set of the first block. 
Then it is clear that the edge set of a general block is a shift of that one,
$$\{ f \in \Pi_1(z_1^k,z_1^{k+1})\cup\Pi_2(z_2^k,z_2^{k+1})\}= z_1^{k-1}+ \mathcal{B}$$ 
By consequence, the variables $(X_i)_{i\leq n}$ are i.i.d and moreover 
\begin{eqnarray*}
\mathbb{E}(X_i)&=& p \big\{ \Psi_0(\theta_S)-\Psi_0(\theta_N)\big\} +q \big\{ \Psi_0(\theta_E)-\Psi_0(\theta_W)\big\}\\
& 
\left\{
\begin{array}{c} =0\\>0\\<0
\end{array}
\right.
&
\begin{array}{cc} 
 {\rm if}  & \theta_N=\theta_S,  \\
 {\rm if}  & \theta_N<\theta_S,  \\
 {\rm if}  & \theta_N>\theta_S.
\end{array}
\end{eqnarray*}
We first consider the case $\theta_N=\theta_S$.
Then $W_k$ is a centered random walk and we can define 
\begin{equation} \label{eq:lnW}
{\mathfrak l}_N= \arg\min_{0 \leq k < N} W_k
\end{equation}
 as in the previous section. 
Before presenting the key lemma,  we introduce the limit law. Denote by $\nu(\cdot | u)$ a regular version of the conditional law
of $(W(a), a\in \mathcal{A} )$ given $\sum_{a \in {\mathcal A}} W(a)=u$. Let $(X_i, i \geq 1)$ be an i.i.d. sequence distributed as
in (\ref{eq:X}), and $S^{\uparrow}$ [resp. $S^{\downarrow}$] associated to $X$ [resp., $-X$] as in (\ref{eq:Sup}), and 
$\hat{S}$ the sequence with $\hat{S}_0=0$ and
  \begin{equation} \nn
 		    \hat{S}_k= \left\{ \begin{array}{ll}
 		    S_k^{\uparrow} & \textrm{if} \ k > 0\\
 		      S_k^{\downarrow} & \textrm{if} \ k<0 \\
 		       \end{array} \right.
 		   \end{equation}
Consider  also, on the same probability space, a random sequence  $(Y_{k,a}: k\in \mathbb{Z}, a\in \mathcal{A})$  such that the vectors
${\mathcal Y}_k=(Y_{k,a}: a\in \mathcal{A})$ are, for $k \in {\mathbb{Z}}$, independent conditionally on $\hat{S}$  with conditional law $\nu(\cdot |  \hat{S}_{k+1} -\hat{S}_k)$.
 \begin{mylem}\label{lemma:p2p}
For fixed $K\in \mathbb{Z}_+ $,
\begin{equation}  \nonumber
\Big[W( \langle z_1^{{\mathfrak l}_N+k},z_2^{{\mathfrak l}_N+k}\rangle+a)-W(\langle z_1^{{\mathfrak l}_N},z_2^{{\mathfrak l}_N}\rangle )\Big]_{|k|\leq K,a\in \mathcal{A}}\stackrel{\mathcal{L}}{\longrightarrow}(\hat{S}_k+Y_{k,a})_{|k|\leq K,a\in \mathcal{A}} .
\end{equation}
\end{mylem}

\textit{Proof.}
 Applying Corollary \ref{col:finiteCV}
to the centered random walk $W_k$ we obtain for fixed $K$
 \begin{equation} \label{eq:dodo}
 \Big[W(\langle z_1^{{\mathfrak l}_N+k},z_2^{{\mathfrak l}_N+k}\rangle)-W(\langle z_1^{{\mathfrak l}_N},z_2^{{\mathfrak l}_N}\rangle)\Big]_{|k|\leq K}\stackrel{\mathcal{L}}{\longrightarrow}(\hat{S}_k)_{|k|\leq K}
 		   \end{equation}
		   with ${\mathfrak l}_N$ from (\ref{eq:lnW}).
 On the other hand, by  independence of the $T_f$'s we know that the vectors 
  \begin{equation} \nonumber
  \Big[W( z_1^{k} +a)-W(\langle z_1^k,z_2^k\rangle)\Big]_{a\in \mathcal{A}} 
  \quad \textmd{are i.i.d. for } k\geq 0. 
   \end{equation}
 Thus, with $W_k=W(\langle z_1^{k},z_2^{k}\rangle)$,  their joint law, conditionally on $(W_k; k \geq 0)$, is $\otimes_k \nu(\cdot | W_{k+1}-W_k)$.
 Then the result follows from (\ref{eq:dodo}).\qed
 \medskip
 
Now we can state the main result of our construction, which reformulates  Theorem \ref{theo:introP2P} in the equilibrium case.

\begin{myprop} \label{th:last} Assume  $\theta_N=\theta_S$.
With the notations of Lemma \ref{lemma:p2p}, let
\begin{equation} \nn
\xi_{k,a}= \frac{\exp\{-\hat{S}_k- Y_{k,a}\}}{1+\sum_{k=-\infty}^{\infty}\sum_{a\in {\mathcal A}} \exp\{-\hat{S}_k- Y_{k,a}\} }
\end{equation}
Then, as $N \to \8$, 
\begin{equation}\label{equa: p2p}
\left(Q_{Np,Nq}^{\omega}(G(\bx)= z_1^{{\mathfrak l}_N+k}+a)\right)_{k\in {\mathbb Z}, a\in \mathcal{A}}\xrightarrow{\mathcal{L}}(\xi_{k,a})_{k\in {\mathbb Z},a\in \mathcal{A}} .
\end{equation}
on the space $(\ell_1(  {\mathbb Z}\times \mathcal{A}), |\cdot|_1)$.
\end{myprop}

\textit{Proof.} 
 We will use the same method as in the section \ref{section: proof-eq} to prove (\ref{equa: p2p}). With Lemma \ref{lemma:p2p}  at hand, we only need here to control the 
 tail of sums as  in (\ref{proof:convergence2b}) and (\ref{proof:CV4b}).  Define
 \begin{eqnarray} 
  \nonumber
(l_{N}^*,a_{N}^*)&=& \arg\min_{(k,a)} \big[W( z_1^{k}+a)\big]\quad  {\rm for \ }  \in \{0,\ldots,N-1\} \times {\mathcal A}, \label{eq:ln*}
\end{eqnarray}
the minimum location of $W$.
Now we consider the process $(k,a) \mapsto 
W( z_1^{k}+a)$ indexed by integer time $t=kN+\ell$ if $a$ is the $\ell$ element in $\mathcal A$, relative to its infimum, i.e., with the shift
$s \mapsto t=s + l_{N}^*\times N+a_{N}^*$. Note that $t \mapsto W( z_1^{k}+a)$ is a sum of independent but not identically distributed random variables, it can be viewed
as a Markov chain,  which is not time-homogeneous but has periodic transitions with period 
equal by the cardinality of $\mathcal A$. Then,  the law
of the post-infimum process 
$$s \mapsto W( z_1^{k}+a)- \min_{m,b} W( z_1^{m}+b), s \geq 0,$$ is also a Markov chain with a lifetime, i.e., a Markov chain killed at a stopping time.
Similar to proposition \ref{theo:finite}, we can prove that the law
of this post-infimum process 
converges as $N \to \8$ to 
a Markov chain on ${\mathbb R}^+$, with non-homogeneous transitions but periodic with period given by the cardinality of $\mathcal A$. The product of $N$ consecutive transition kernels  
coincides with the one of $S^{\uparrow}$, it  is homogeneous. Similar to Theorem \ref{th:ritter}, we conclude that the post-infimum process grows algebraically: with probability arbitrarily close to 1, we have for some positive $\delta$  and all large $N$,
\begin{equation}
 \nonumber
W( z_1^{l_{N}^*+k}+a)- \min_{m,b} W( z_1^{m}+b) \geq \delta k^\eta \;, \quad k \in \{1,\ldots, N-1-l_{N}^*\}.
\end{equation}
Then, it is plain to check that ${\mathfrak l}_N-l_{N}^*= {\mathcal O}(1)$ in probability using that the former minimizes $W_k=W(\langle z_1^k, z_2^k\rangle)$, and we derive (\ref{proof:convergence2b}) and (\ref{proof:CV4b}) as well. The rest of the proof follows from
similar arguments  to those of  Theorem \ref{theo:main} in the case of  $\theta=\mu/2$ and from  Lemma \ref{lemma:p2p}. \qed
\medskip

\textit{ 
Proof of  Theorem \ref{theo:introP2P}:} In the equilibrium case $\theta_N=\theta_S$, the above proposition \ref{th:last} yields the conclusion
by taking
$$
m_N =  (z_1^{{\mathfrak l}_N}+z_2^{{\mathfrak l}_N}) \cdot  (q{\bf e_1}-p{\bf e_2}) ,
$$
so that any path $\bx$ through  the edge $\langle z_1^{{\mathfrak l}_N},z_2^{{\mathfrak l}_N} \rangle$ has $F(\bx)=m_M$.

 Consider now the opposite case where, by symmetry,  we may assume $\theta_N<\theta_S$ without loss of generality. Then, the walk $(W_k)_{k \geq 0}$ has a 
 global minimum. 
Interpolating $(W_k)_{k \geq 0}$ with independent pieces with law $\nu$ defined as above, we construct a process  $\hat \xi$. Repeating the 
arguments of section \ref{sec:4.2},  we check that it is the desired limit. \qed

\medskip

Similar to that of Theorem \ref{theo:arcsin},
the proof of Theorem \ref{theo:arcsin2} is straightforward and left to the reader.

\medskip
{\em Acknowledgements:} We thank the referee for 
a careful reading, many useful comments
and numerous suggestions to improve the paper.

{ 
\footnotesize
 \bibliographystyle{plain}
\bibliography{gamma}

\begin{thebibliography}{10}

\bibitem{Bert91b}
J.~Bertoin.
\newblock D\'ecomposition du mouvement brownien avec d\'erive en un minimum
  local par juxtaposition de ses excursions positives et n\'egatives.
\newblock In {\em S\'eminaire de {P}robabilit\'es, {XXV}}, volume 1485 of {\em
  Lecture Notes in Math.}, pages 330--344. Springer, Berlin, 1991.

\bibitem{Bert91a}
J.~Bertoin.
\newblock Sur la d\'ecomposition de la trajectoire d'un processus de {L}\'evy
  spectralement positif en son infimum.
\newblock {\em Ann. Inst. H. Poincar\'e Probab. Statist.}, 27(4):537--547,
  1991.

\bibitem{Bert93}
J.~Bertoin.
\newblock Splitting at the infimum and excursions in half-lines for random
  walks and {L}\'evy processes.
\newblock {\em Stochastic Process. Appl.}, 47(1):17--35, 1993.

\bibitem{BeDo94}
J.~Bertoin and R.~A. Doney.
\newblock On conditioning a random walk to stay nonnegative.
\newblock {\em Ann. Probab.}, 22(4):2152--2167, 1994.

\bibitem{Bolt89}
E.~Bolthausen.
\newblock A note on the diffusion of directed polymers in a random environment.
\newblock {\em Comm. Math. Phys.}, 123(4):529--534, 1989.

\bibitem{BoroCorwReme13}
A.~Borodin, I.~Corwin, and D.~Remenik.
\newblock Log-gamma polymer free energy fluctuations via a {F}redholm
  determinant identity.
\newblock {\em Comm. Math. Phys.}, 324(1):215--232, 2013.

\bibitem{CaHu02}
P.~Carmona and Y.~Hu.
\newblock On the partition function of a directed polymer in a {G}aussian
  random environment.
\newblock {\em Probab. Theory Related Fields}, 124(3):431--457, 2002.

\bibitem{CaMo94}
R.~Carmona and S.~Molchanov.
\newblock Parabolic {A}nderson problem and intermittency.
\newblock {\em Mem. Amer. Math. Soc.}, 108(518), 1994.

\bibitem{CoCr13}
F.~Comets and M.~Cranston.
\newblock Overlaps and pathwise localization in the {A}nderson polymer model.
\newblock {\em Stochastic Process. Appl.}, 123(6):2446--2471, 2013.

\bibitem{CoShYo03}
F.~Comets, T.~Shiga, and N.~Yoshida.
\newblock Directed polymers in a random environment: path localization and
  strong disorder.
\newblock {\em Bernoulli}, 9(4):705--723, 2003.

\bibitem{CoYo06}
F.~Comets and N.~Yoshida.
\newblock Directed polymers in random environment are diffusive at weak
  disorder.
\newblock {\em Ann. Probab.}, 34(5):1746--1770, 2006.

\bibitem{CoYo13}
F.~Comets and N.~Yoshida.
\newblock Localization transition for polymers in {P}oissonian medium.
\newblock {\em Comm. Math. Phys.}, 323(1):417--447, 2013.

\bibitem{CorwinKPZ12}
I.~Corwin.
\newblock The {K}ardar-{P}arisi-{Z}hang equation and universality class.
\newblock {\em Random Matrices Theory Appl.}, 1(1):1130001, 76, 2012.

\bibitem{CoOcSeppZyg14}
I.~Corwin, N.~O'Connell, T.~Sepp{\"a}l{\"a}inen, and N.~Zygouras.
\newblock Tropical combinatorics and {W}hittaker functions.
\newblock {\em Duke Math. J.}, 163(3):513--563, 2014.

\bibitem{DamronHanson14}
M.~Damron and J.~Hanson.
\newblock Busemann functions and infinite geodesics in two-dimensional
  first-passage percolation.
\newblock {\em Comm. Math. Phys.}, 325(3):917--963, 2014.

\bibitem{Doney89}
R.~A. Doney.
\newblock Last exit times for random walks.
\newblock {\em Stochastic Process. Appl.}, 31(2):321--331, 1989.

\bibitem{Feller}
W.~Feller.
\newblock {\em An introduction to probability theory and its applications.
  {V}ol. {II}.}
\newblock Second edition. John Wiley \& Sons, Inc., New York-London-Sydney,
  1971.

\bibitem{GaPeSh10}
N.~Gantert, Y.~Peres, and Z.~Shi.
\newblock The infinite valley for a recurrent random walk in random
  environment.
\newblock {\em Ann. Inst. Henri Poincar\'e Probab. Stat.}, 46(2):525--536,
  2010.

\bibitem{GeRASe14}
N.~Georgiou, F.~Rassoul-Agha, and T.~Sepp{\"a}l{\"a}inen.
\newblock Stationary cocycles for the corner growth model.
\newblock {\em preprint, arXiv: 1404.7786}, 2014.

\bibitem{GeRASeYi14}
N.~Georgiou, F.~Rassoul-Agha, T.~Sepp{\"a}l{\"a}inen, and A.~Yilmaz.
\newblock Ratios of partition functions for the log-gamma polymer.
\newblock {\em preprint, arXiv:1303.1229}, 2013.

\bibitem{GeorSepp13}
N.~Georgiou and T.~Sepp{\"a}l{\"a}inen.
\newblock Large deviation rate functions for the partition function in a
  log-gamma distributed random potential.
\newblock {\em Ann. Probab.}, 41(6):4248--4286, 2013.

\bibitem{Golo84}
A.~O. Golosov.
\newblock Localization of random walks in one-dimensional random environments.
\newblock {\em Comm. Math. Phys.}, 92(4):491--506, 1984.

\bibitem{Hairer13}
M.~Hairer.
\newblock Solving the {KPZ} equation.
\newblock {\em Ann. of Math. (2)}, 178(2):559--664, 2013.

\bibitem{HaKeKy03}
B.~M. Hambly, G.~Kersting, and A.~E. Kyprianou.
\newblock Law of the iterated logarithm for oscillating random walks
  conditioned to stay non-negative.
\newblock {\em Stochastic Process. Appl.}, 108(2):327--343, 2003.

\bibitem{HuHe85}
D.A. Huse and C.L. Henley.
\newblock Pinning and roughening of domain wall in {I}sing systems due to
  random impurities.
\newblock {\em Phys. Rev. Lett.}, 54:2708--2711, 1985.

\bibitem{ImSp88}
J.~Z. Imbrie and T.~Spencer.
\newblock Diffusion of directed polymers in a random environment.
\newblock {\em J. Statist. Phys.}, 52(3-4):609--626, 1988.

\bibitem{KaSh98}
I.~Karatzas and S.~Shreve.
\newblock {\em Methods of mathematical finance}, volume~39 of {\em Applications
  of Mathematics (New York)}.
\newblock Springer-Verlag, New York, 1998.

\bibitem{KeMe04}
G.~Kersting and K.~Memi{\c{s}}o{\v{g}}lu.
\newblock Path decompositions for {M}arkov chains.
\newblock {\em Ann. Probab.}, 32(2):1370--1390, 2004.

\bibitem{Levy39}
P.~L{\'e}vy.
\newblock Sur certains processus stochastiques homog\`enes.
\newblock {\em Compositio Math.}, 7:283--339, 1939.

\bibitem{MoOc07}
J.~Moriarty and N.~O'Connell.
\newblock On the free energy of a directed polymer in a {B}rownian environment.
\newblock {\em Markov Process. Related Fields}, 13(2):251--266, 2007.

\bibitem{Newman97}
C.~Newman.
\newblock {\em Topics in disordered systems}.
\newblock Lectures in Mathematics ETH Z\"urich. Birkh\"auser Verlag, Basel,
  1997.

\bibitem{OcYor01}
N.~O'Connell and M.~Yor.
\newblock Brownian analogues of {B}urke's theorem.
\newblock {\em Stochastic Process. Appl.}, 96(2):285--304, 2001.

\bibitem{QuastellKPZ10}
J.~Quastel.
\newblock Weakly asymmetric exclusion and {KPZ}.
\newblock In {\em Proceedings of the {I}nternational {C}ongress of
  {M}athematicians. {V}olume {IV}}, pages 2310--2324. Hindustan Book Agency,
  New Delhi, 2010.

\bibitem{Ritter81}
G.~A. Ritter.
\newblock Growth of random walks conditioned to stay positive.
\newblock {\em Ann. Probab.}, 9(4):699--704, 1981.

\bibitem{Sepp12}
T.~Sepp{\"a}l{\"a}inen.
\newblock Scaling for a one-dimensional directed polymer with boundary
  conditions.
\newblock {\em Ann. Probab.}, 40(1):19--73, 2012.

\bibitem{Sinai82}
Ya.~G. Sina{\u\i}.
\newblock The limit behavior of a one-dimensional random walk in a random
  environment.
\newblock {\em Teor. Veroyatnost. i Primenen.}, 27(2):247--258, 1982.

\bibitem{Sina95}
Ya.~G. Sina{\u\i}.
\newblock A remark concerning random walks with random potentials.
\newblock {\em Fund. Math.}, 147(2):173--180, 1995.

\bibitem{Tana89}
H.~Tanaka.
\newblock Time reversal of random walks in one-dimension.
\newblock {\em Tokyo J. Math.}, 12(1):159--174, 1989.

\bibitem{Varg06}
V.~Vargas.
\newblock A local limit theorem for directed polymers in random media: the
  continuous and the discrete case.
\newblock {\em Ann. Inst. H. Poincar\'e Probab. Statist.}, 42(5):521--534,
  2006.

\bibitem{Vargas07}
V.~Vargas.
\newblock Strong localization and macroscopic atoms for directed polymers.
\newblock {\em Probab. Theory Related Fields}, 138(3-4):391--410, 2007.

\bibitem{Will74}
D.~Williams.
\newblock Path decomposition and continuity of local time for one-dimensional
  diffusions. {I}.
\newblock {\em Proc. London Math. Soc. (3)}, 28:738--768, 1974.

\end{thebibliography}
}

\end{document}